\documentclass{article}
\usepackage[utf8]{inputenc}
\usepackage[T1]{fontenc}
\usepackage{amsthm}
\usepackage{amssymb}
\usepackage{mathrsfs}
\usepackage{amsmath}
\usepackage{amsfonts,amssymb}
\usepackage{mathtools}
\usepackage{fullpage}
\usepackage{hyperref}
\usepackage{multirow}
\usepackage{relsize}
\usepackage{lmodern}
\usepackage{import}
\usepackage{authblk}

\usepackage[usenames, dvipsnames]{xcolor}
\usepackage{tikz}
\usepackage{tikz-3dplot}
\usetikzlibrary{arrows}
\usetikzlibrary{decorations}
\usetikzlibrary{matrix}
\usetikzlibrary{calc}
\newcommand{\midarrow}{\tikz \draw[-triangle 90] (0,0) -- +(.1,0);}
\usepackage{ifthen}

\definecolor{red}{rgb}{1,0.25,0.2}
\definecolor{blue}{rgb}{0,0,0.5}
\definecolor{green}{rgb}{0,0.8,0}

\newtheorem{Theo}{Theorem}
\newtheorem*{Theo*}{Theorem}
\newtheorem{Prop}{Proposition}

\newtheorem{Le}{Lemma}
\newtheorem{Rk}{Remark}
\theoremstyle{definition}
\newtheorem{Def}{Definition}
\newtheorem*{Defs*}{Definitions}
\newtheorem*{Ex}{Example}

\newcommand*{\Scale}[2][4]{\scalebox{#1}{$#2$}}%

\newcommand{\R}{\mathbb{R}}
\newcommand{\C}{\mathbb{C}}
\renewcommand{\Pr}{\mathbb{P}}
\newcommand{\N}{\mathbb{N}}
\newcommand{\Z}{\mathbb{Z}}
\newcommand{\T}{\mathbb{T}}
\newcommand{\ZZ}{\mathcal{Z}}
\newcommand{\YY}{\mathcal{Y}}
\newcommand{\F}{\mathcal{F}}

\newcommand{\loc}{\mathrm{loc}}

\title{The free-fermionic $C^{(1)}_2$ loop model, double dimers and
  Kashaev's recurrence}
\author{Paul Melotti\thanks{Laboratoire de
    Probabilités, Statistique et Modélisation, Sorbonne Université,
    campus Pierre et Marie Curie, 4 place Jussieu, F-75005
    Paris. Email: \texttt{paul.melotti@upmc.fr}}}
\date{\today}

\begin{document}

\maketitle

\begin{abstract} We study a two-color loop model known as the
  $C^{(1)}_2$ loop model. We define a free-fermionic regime for this
  model, and show that under this assumption it can be transformed
  into a double dimer model. We then compute its free energy on
  periodic planar graphs. We also study the star-triangle relation or
  Yang-Baxter equations of this model, and show that after a proper
  parametrization they can be summed up into a single relation known
  as Kashaev's relation. This is enough to identify the solution of
  Kashaev's relation as the partition function of a $C^{(1)}_2$ loop
  model with some boundary conditions, thus solving an open question
  of Kenyon and Pemantle \cite{KenyonPemantle} about the combinatorics
  of Kashaev's relation.
\end{abstract}


\section{Introduction}

In 1996 Kashaev introduced a way to rewrite the star triangle
transformation of the Ising model \cite{Kashaev}. Specifically, let us
take a planar graph $G=(V,E)$ with usual coupling constants for the
Ising model $(J_e)_{e\in E}$ on the edges. Let us suppose that there
is a set of variables $g$ on the vertices and faces of $G$ such that
\begin{equation}
  \label{eq:kashis} \sinh^2(J_e) = \frac{g_xg_y}{g_ug_v}
\end{equation}
where $x,y$ are the endpoints of $e$ and $u,v$ are the
faces adjacent to $e$.  Then the star-triangle relation, or local
Yang-Baxter equation, is equivalent to the variables $g$ satisfying a
single polynomial relation:

\begin{equation*}
  \raisebox{-1.2cm}{\input{tikz/st0.tex}}
  \hspace*{-2.5cm}
  \begin{split} & g^2 g_{123}^2 + g_{1}^2 g_{23}^2 + g_{2}^2 g_{13}^2
    + g_{3}^2 g_{12}^2 \\ - & 2 g_{2}g_{3}g_{13}g_{12} - 2
    g_{1}g_{3}g_{23}g_{12} -2 g_{1}g_{2}g_{23}g_{13} \\ - & 2 g g_{123}
    (g_{1} g_{23} + g_{2} g_{13} + g_{3} g_{12} ) \\ - & 4 g g_{12} g_{23}
    g_{13} - 4 g_{123} g_1 g_2 g_3 \\ =& 0.
  \end{split}
\end{equation*}

This relation\footnote{Kashaev's initial equation contained a $+4$
  instead of a $-4$ coefficient for the last terms, but one can easily
  get from one to another, for instance by multiplying $g$ by $-1$ at
  a vertex of the cube and its three neighbors.}, known as Kashaev's
relation, has sparked some interest from the point of view of
\textit{spatial recurrences}. It can be embedded in $\Z^3$ by taking
$x\in \Z^3$ and denoting $g = g_x$, $g_i = g_{x+e_i}$,
$g_{ij}=g_{x+e_i+e_j}$, etc. Then by choosing the greatest root of a
degree 2 polynomial we get \cite{KenyonPemantle}:

\begin{equation}
  \label{eq:krec} g_{123} = \frac{2g_1 g_2 g_3 + g(g_1 g_{23} + g_2
    g_{13} + g_3 g_{12}) +2XYZ}{g^2},
\end{equation}
where $X = \sqrt{g g_{23} + g_2 g_3}$,
$Y = \sqrt{g g_{13} + g_1 g_3}$, $Z = \sqrt{g g_{12} + g_1 g_2}$.

This transformation \eqref{eq:krec} is called \textit{Kashaev's
  recurrence}. It can be iterated to define $g$ on further vertices of
$\Z^3$, provided we had a sufficiently large set of initial
conditions. A remarkable fact is that it exhibits a
\textit{Laurentness} phenomenon: the solution of the recurrence at any
point is always a Laurent polynomial in the initial variables. This
fact is related to cluster algebras \cite{FZ1,FZ2}, but it also hints
at a possible hidden object represented by the solution.

Let us quickly review the current state of spatial recurrences: Speyer
related the solution of the octahedron recurrence (which can be traced
back to Dodgson \cite{Dodgson}) to the partition function of a dimer
model \cite{Speyer_octa}; then Carroll and Speyer showed that the cube
recurrence (proposed by Propp \cite{Propp}) corresponds to cube groves
\cite{CarrollSpeyer}; more recently Kenyon and Pemantle studied a
generalization of Kashaev's relation, known as the hexahedron
recurrence, and identified its solution with a double dimer model
\cite{KenyonPemantle}. Unfortunately, when specialized to Kashaev's
recurrence, their model does not provide a one-to-one correspondence
between configurations and monomials of the Laurent polynomial.  In
this paper we provide a model that does give a one-to-one
correspondence, known in the physics literature as the $C^{(1)}_2$
loop model.

\medskip

The $C^{(1)}_2$ loop model was introduced by Warnaar and Nienhuis in
\cite{WarnaarNienhuis}, among other models, as a loop model naturally
associated to formal solutions of the Yang-Baxter equation
\cite{Bazhanov,Jimbo}. It was also considered, for different reasons,
by Jacobsen and Kondev in \cite{JacobsenKondev} as a generalization of
the eight-vertex model, and they conjecture a phase diagram for this
model.

\begin{figure}[!h] \centering \def\svgwidth{3.5cm}
  \import{./}{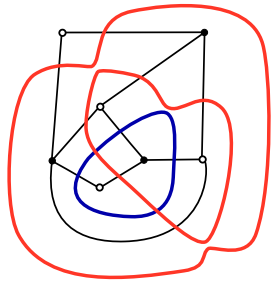}
  \caption{A quadrangulation $G$ (in black), with a $C^{(1)}_2$ loop
    configuration.}
  \label{fig:qloop}
\end{figure}

It is a dense, two-color loop model where same-colored loops cannot
intersect. Let us detail this terminology; see also
Figure~\ref{fig:qloop} for an example.

\begin{figure}[!h] \centering \input{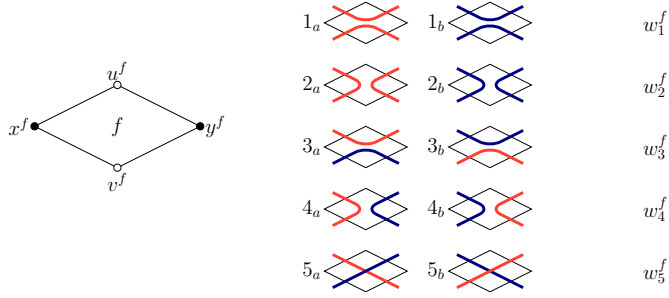}
  \caption{A face $f$, the 10 local configurations $1_a, \dots, 5_b$
    at $f$ and their local weight $w^f_i$.}
  \label{fig:locconf}
\end{figure}

\begin{itemize}
\item \textit{loop model} means that the configurations are unions of
  simple curves. In our case, these curves use edges of the dual graph
  of a bipartite quadrangulation $G$, and are able to turn inside
  faces of $G$.
\item \textit{two-color} means that each loop is either red or blue.
\item \textit{dense} means that every dual edge of $G$ belongs to a
  loop.
\item finally, \textit{same-colored loops cannot intersect} means that
  the only allowed crossings are between red and blue; see
  Figure~\ref{fig:locconf} for all possible local configurations at a
  face of $G$.
\end{itemize}

Equivalently, a $C^{(1)}_2$ loop configuration can be seen as a gluing
of quadrangles appearing in Figure~\ref{fig:locconf}, where only
same-colored edges can be glued together.

The model is equipped with weights: let $n$ be a positive parameter
called \textit{fugacity}, then the weight of a $C^{(1)}_2$ loop
configuration $\sigma$ is
\begin{equation*}
  w(\sigma) = n^{\# \mathrm{loops \ in \ }\sigma} \ 
  \prod_{f\in F} w^f_i,
\end{equation*}
where $F$ is the set of faces of the quadrangulation and $w^f_i$ is
the local weight corresponding to $\sigma|_{f}$.

In \cite{WarnaarNienhuis} the authors give a two-parameter family of
fugacity and weights such that the model is \textit{integrable}, which
in their setting means that it satisfies a form of the star-triangle
relation:

\begin{equation}
  \label{eq:pdswn}
  \begin{cases} n = -2 \cos 2\lambda \\
    w_1 = \frac{\sin(\lambda-u)
      \sin(3\lambda - u)}{\sin \lambda \sin 3\lambda}\\
    w_2 = - \frac{\sin u
      \sin(2\lambda - u)}{\sin \lambda \sin 3\lambda}\\
    w_3 = \frac{
      \sin(3\lambda - u)}{\sin 3\lambda}\\
    w_4 = - \frac{\sin u}{\sin
      3\lambda}\\
    w_5 = \frac{\sin u \sin(3\lambda - u)}{\sin \lambda \sin
      3\lambda}
  \end{cases}
\end{equation}

In \cite{IkhlefCardy}, Ikhlef and Cardy define a fermionic observable
$F_s$ for this model and show that imposing a form of discrete
holomorphicity on $F_s$ yields the same integrable weights as in
\cite{WarnaarNienhuis}; this approach was extended to the case of
non-trivial boundary conditions in \cite{deGierLeeRasmussen}.

\medskip

In this paper we only deal with the $n=2$ case. We introduce the
\textit{free-fermionic} relations (see Section~\ref{sec:ff} for more
details on this terminology):
\begin{equation*}
  \begin{split}
    w^f_1 w^f_4 & = w^f_3 w^f_5,\\
    w^f_2 w^f_3 & = w^f_4 w^f_5, \\
    w^f_5 (w^f_1 + w^f_2) & = w^f_3 w^f_4.
  \end{split}
\end{equation*}

For instance, the integrable weights \eqref{eq:pdswn} at $n=2$
(\textit{i.e.} $\lambda = \pm \frac{\pi}{2}$) satisfy the
free-fermionic relations.

When these relations are satisfied, we show that the $C^{(1)}_2$ loop
model can be transformed into a double dimer model. We prove the
following; for a precise statement, see Theorem~\ref{theo:c12dim} of
Section~\ref{sec:loopsdimers}.
\begin{Theo*} 
  For any free-fermionic $C^{(1)}_2$ loop model, there is
  a bipartite decorated graph equipped with a dimer model with weights
  $\mu$, and there are constants $(\lambda_f)_{f\in F}$, such that the
  partition function $\ZZ_{\text{loop}}^G$ is equal to the square of
  the partition function $\ZZ_{\text{dim}}\left(\mu\right)$ of this
  dimer model, up to multiplicative factors:
  \begin{equation*} \ZZ_{\text{loop}}^G = \left(\prod_{f\in F}
      \lambda_f\right)
    \left(\ZZ^G_{\text{dim}}\left(\mu\right)\right)^2.
  \end{equation*}
\end{Theo*}

An application of this result is the computation of the free energy of
any free-fermionic $C^{(1)}_2$ loop model on a periodic planar
quadrangulation; see Section~\ref{sec:freeen}.

Then we define a parametrized free-fermionic $C^{(1)}_2$ loop model,
analogous to Kashaev's parametrized Ising model \eqref{eq:kashis}. Let
us suppose that there is a set of variables $(g_v)_{v\in V}$ on the
vertices of $G$ such that

\begin{equation}
  \label{eq:pdsg0}
  \begin{cases} w^f_1 = g_xg_y \\
    w^f_2 = g_u g_v \\
    w^f_3 = \sqrt{g_xg_y} \sqrt{g_xg_y + g_ug_v} \\
    w^f_4 = \sqrt{g_vg_u} \sqrt{g_xg_y + g_ug_v}\\
    w^f_5 = \sqrt{g_xg_yg_ug_v}.
  \end{cases}
\end{equation}
The existence of such a parametrization is discussed in
Appendix~\ref{app:fullc12}. In particular, we show that it always
exists for a free-fermionic model on a lozenge graph.

In this regime, we show that the Yang-Baxter equations associated to
the model (corresponding to a move called \textit{cube flip}, similar
to the star/triangle move) are equivalent to $g$ satisfying Kashaev's
recurrence \eqref{eq:krec}. See Theorem~\ref{theo:cuberem} in
Section~\ref{sec:st}.

\begin{figure}[!h] 
  \centering 
  \input{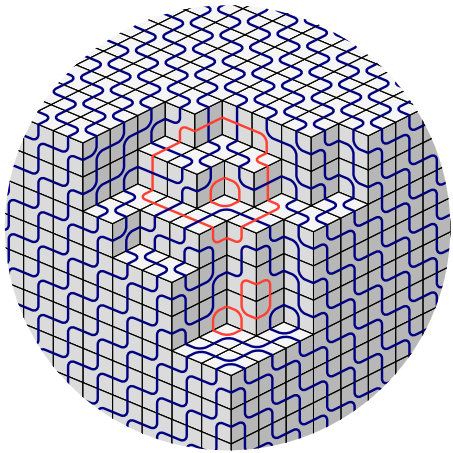}
  \caption{A taut $C^{(1)}_2$ loop configuration.}
  \label{fig:ex}
\end{figure}

Finally, we get to the solution of Kashaev's recurrence. See
Theorems~\ref{theo:princ2} and \ref{theo:unic} of
Section~\ref{sec:taut} for a precise statement of the following.
\begin{Theo*} 
  For any solvable initial condition 
  $(g_i)_{i\in I}$ on
  $I \subset \mathbb{Z}_{-} ^3$, the solution of Kashaev's recurrence at
  the origin is
  \begin{equation*} 
    g_{0,0,0} = \sum_{\sigma} \left( 2^{\# \mathrm{loops \ in \ }\sigma} \ 
    \prod_{f} w^f_i \ \prod_{i\in I} g_i^{-2} \right)
  \end{equation*}
  where the sum is over \textit{taut} $C^{(1)}_2$ configurations
  $\sigma$, and the local weights $w^f_i$ are given by
  \eqref{eq:pdsg0}.

  Moreover, there is a one-to-one correspondence between such loop
  configurations and monomials of $g_{0,0,0}$ as a function of the
  variables $(g_i)_{i\in I}$.
\end{Theo*}
The \textit{taut} configurations are those satisfying some boundary
and connectivity conditions; an example is displayed in
Figure~\ref{fig:ex}.

The paper is organized as follows:

\begin{itemize}
\item In Section~\ref{sec:loopsdimers} we define the models, and show
  how to get from free-fermionic $C^{(1)}_2$ loops to dimers. We also
  compute the free energy on periodic planar quadrangulations.
\item In Section~\ref{sec:st} we show that for the parametrization
  \eqref{eq:pdsg0} of the weights, the Yang-Baxter equations of loops
  are equivalent to Kashaev's recurrence.
\item In Section~\ref{sec:taut} we define taut configurations, and
  prove that the solution of Kashaev's recurrence is the partition
  function of these configurations. We compute some limit shapes of
  the model by a now standard technique
  \cite{PetersenSpeyer,diFrancescoSotoGarrido,KenyonPemantle,George}. We
  also show that in characteristic $2$ the model reduces to the cube
  groves of \cite{CarrollSpeyer}.
\end{itemize}

\subsection*{Acknowledgments}
I am very grateful to my Ph.D. advisors Cédric Boutillier and Béatrice
de Tilière for the motivations and questions related to the present
work, and for their constant help during the writing of this paper. I
also thank an anonymous referee for pointing at a mistake in the
definition of the free-fermionic regime.

\section{Free-fermionic $C^{(1)}_2$ loops and double dimers}
\label{sec:loopsdimers}

\subsection{The $C^{(1)}_2$ loop model on a quadrangulation}

\begin{Def}
  \label{def:quadrangulation} Let $\mathcal{S}$ be a connected
  orientable compact surface without boundary. A
  \textit{quadrangulation} of $\mathcal{S}$ is a finite connected
  simple graph $G=(V,E)$ embedded in $\mathcal{S}$ so that edges do
  not intersect, and so that the \textit{faces} of $G$ (the connected
  components of the complement of the embedding) are homeomorphic to
  disks and have degree $4$. We denote by $F$ the set of faces.
\end{Def}

\begin{Def} Let $G$ be a bipartite quadrangulation of
  $\mathcal{S}$. For every face $f \in F$, we fix names for the
  vertices of the boundary of $f$, in clockwise order, as
  $x^f,u^f,y^f,v^f$, with $x^f,y^f$ black vertices and $u^f,v^f$ white
  vertices like in Figure~\ref{fig:locconf}. Notice that a vertex will
  have several names, corresponding to all the faces adjacent to it;
  we only use these labels to make the $10$ different configurations
  of Figure~\ref{fig:locconf} well defined. When there is no
  ambiguity, we will also drop the superscript $f$.

  A $C^{(1)}_2$ loop configuration $\sigma$ on $G$ is the data, for
  every $f\in F$, of an index
  $i^f_k \in \{1_a, 1_b, \dots, 5_a, 5_b\}$ (we think of it as
  $i\in\{1,\dots,5\}$ and $k\in\{a,b\}$) representing the local
  configuration $\sigma|_f$, such that glued edges are the same color.

  When there is no ambiguity on the face involved, we will often drop
  the superscript $f$ in $i^f_k$. Let us denote $\sigma$ the set
  $\sigma=\left\{ (f,i_k) \mid f\in F \right\}$.
\end{Def}

Simply stated, a $C^{(1)}_2$ loop configuration on $G$ is an
edge-covering set of red or blue loops on the dual graph $G^*$ (which
we represent inside a face by turning or crossing when necessary),
such that same-colored loops cannot cross.

We equip the faces with a set of positive weights
$W = (w^f_i)_{f\in F, i\in \{1,\dots,5\}}$. For a loop configuration
$\sigma$ we let $N_{\sigma}$ be the number of loops in $\sigma$. The
weight of $\sigma$ is defined as

\begin{equation}
  \label{eq:weight2} w_{\text{loop}}^G(\sigma) = 2^{N_{\sigma}}
  \prod_{(f,i_k)\in \sigma}w^f_i.
\end{equation}

The \textit{partition function} of the model is the weighted sum of
loop configurations:
\begin{equation*} \ZZ_{\text{loop}}^G \left(W\right) = \sum_{\sigma}
  w_{\text{loop}}^G(\sigma).\end{equation*}

\subsection{Dimer model on the quad-graph} Let $G^*$ be the dual graph
of $G$. We consider a decorated graph, denoted $G^Q$, constructed by
expanding every vertex of $G^*$ (which has degree 4) into a small
quadrangle called a \textit{city}\footnote{This terminology is related
  to the \textit{urban renewal} transformation of the dimer model
  \cite{KenyonProppWilson}.}.  Cities are connected by edges called
\textit{roads}. Let $E^Q$ be the edges of $G^Q$.

\begin{figure}[!h]
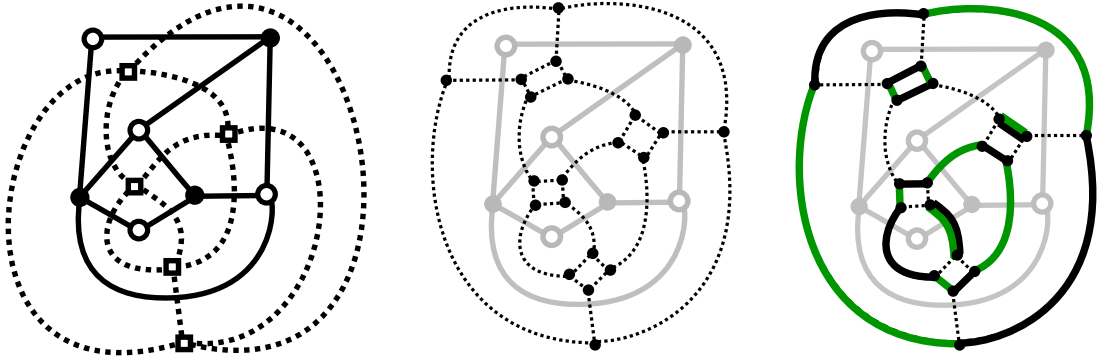
 \centering \def\svgwidth{4.8cm}
  \import{./}{fig/dual.tex} \hspace*{0.5cm} \def\svgwidth{4cm}
  \import{./}{fig/quad.tex} \hspace*{0.5cm} \def\svgwidth{4cm}
  \import{./}{fig/dim.tex}
  \caption{A quadrangulation $G$ of the sphere with its dual $G^*$
    (dotted); the decorated graph $G^Q$ (dotted); a double dimer model
    $(m_A,m_B)$ on $G^Q$.}
  \label{fig:quad}
\end{figure}

A dimer configuration on $G^Q$ is a subset $m\subset E^Q$ such that
every vertex of $G^Q$ belongs to exactly one edge of $m$. Dimer models
on $G^Q$ have appeared several times
\cite{WuLin,deTiliere,Dubedat,BoutillierDeTiliere} in the study of the
6V model and of the Ising model.

Let $\mu=(\mu_e)_{e\in E^Q}$ be a set of positive real weights on the
edges of $G^Q$. The weight of a dimer configuration is then:
\begin{equation*}
  w^G_{\text{dim}}(m) = \prod_{e\in m}\mu_e.
\end{equation*}

We similarly define the partition function for dimers:
\begin{equation*}
  \ZZ^G_{\text{dim}}\left( \mu \right) = \sum_{m} w^G_{\text{dim}}(m),
\end{equation*}
where the sum is over all dimer configurations of $G^Q$.

The aim of this Section is to provide a direct link between the
$C^{(1)}_2$ loop model on $G$ in a certain regime and a couple of
independent dimer models on $G^Q$.

\subsection{Free-fermion regime}
\label{sec:ff}

Let us make the following assumptions on the $C^{(1)}_2$ loops weights
(implicitly evaluated at a face $f\in F$), which we call the
\textit{free-fermionic relations}:
\begin{equation}
  \label{eq:ff}
  \begin{split} w_1 w_4 & = w_3 w_5,\\
    w_2 w_3 & = w_4 w_5, \\
    w_5(w_1+w_2) & = w_3 w_4.
  \end{split}
\end{equation}
Let us make a few remarks on this terminology. We will
see that relations \eqref{eq:ff} are sufficient to transform the
$C^{(1)}_2$ loop model into dimers. This idea has been used several
times in statistical mechanics to get exact solutions for various
models, such as the Ising model \cite{Kasteleyn:ising} and various
vertex models \cite{FanWu,FanWu2,Assis}. In the physics literature,
this technique is sometimes called the ``Pfaffian method'', since the
dimer model's partition function corresponds to Pfaffians
\cite{Kasteleyn:ising}. An alternative representation of Pfaffians is
to use Grassman integrals, which are integrals of anti-commuting
variables; see for instance \cite{DiFrancescoEtAl}, chapter 2.B. These
anti-commuting variables are interpreted physically as a system
non-interacting fermions \cite{Hurst}. This is why any regime for
which there is a transformation to dimers is often called
\textit{free-fermionic}.

\begin{Le}
  Let $w_1,w_2,w_3,w_4,w_5 \in (0,\infty)$ be five positive real
  numbers. Then they satisfy \eqref{eq:ff} \textit{iff} there exists a
  unique triplet $\lambda,a,b \in (0,\infty)$ such that
  $a^2 + b^2 = 1$ and
  \begin{equation}
    \label{eq:pdsff}
    \begin{cases}
      w_1 = \lambda a^2\\
      w_2 = \lambda b^2\\
      w_3 = \lambda a\\
      w_4 = \lambda b\\
      w_5 = \lambda ab.
    \end{cases}
  \end{equation}
\end{Le}

\begin{proof} Given a set of weights $w_1,w_2,w_3,w_4,w_5 \in
  (0,\infty)$ that satisfy \eqref{eq:ff}, then there is only one
  candidate for $\lambda,a,b$:
  \begin{equation*}
    \begin{split} \lambda & = w_1 + w_2, \\
      a & = \frac{w_3}{\lambda},\\
      b & = \frac{w_4}{\lambda}.
    \end{split}
  \end{equation*}
  Then the third equation of \eqref{eq:ff} simplifies into
  \begin{equation*}
    w_5 = \lambda a b
  \end{equation*}
  and the first two equations become, respectively
  \begin{equation*}
    \begin{split}
      w_1 & = \lambda a^2,\\
      w_2 & = \lambda b^2
    \end{split}
  \end{equation*}
  and since $\lambda = w_1+w_2$, we also have $a^2 + b^2 = 1$ so that
  the parameterization \eqref{eq:pdsff} is correct. Reciprocally, it
  is easy to check that \eqref{eq:pdsff} implies \eqref{eq:ff}.
\end{proof}

\begin{figure}[h]
  \centering \input{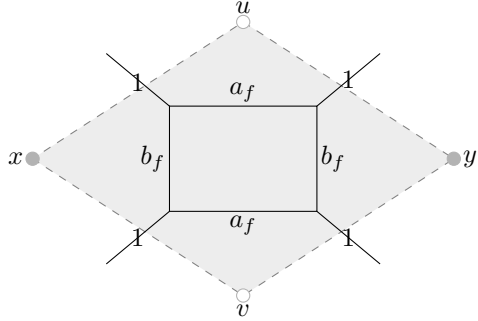}
  \caption{Edge weights for the dimer model on $G^Q$ at a face $f\in
    F$.}
  \label{fig:dimweightgen}
\end{figure}

\begin{Theo}
  \label{theo:c12dim} Let us consider a symmetric $C^{(1)}_2$ loop
  model such that the free-fermionic relations \eqref{eq:ff} are
  satisfied at each face $f\in F$. Let $\lambda_f,a_f,b_f$ be the
  corresponding parameters in representation \eqref{eq:pdsff}.
  
  Let us consider the dimer model on $G^Q$ with weights
  $\mu=(\mu_e)_{e\in E^Q}$ given by Figure~\ref{fig:dimweightgen}.
  
  Then
  \begin{equation}
    \label{eq:zloopdim} \ZZ_{\text{loop}}^G \left( W \right) =
    \left(\prod_{f\in F} \lambda_f\right)
    \left(\ZZ^G_{\text{dim}}\left(\mu\right)\right)^2.
  \end{equation}
\end{Theo}

\begin{proof}

  \subsubsection*{Step 1: double dimers and face weights} Let us
  consider a couple of independent dimer configurations $m_A$ and
  $m_B$. Dimers $A$ will be colored in black and dimers $B$ in
  green. Clearly,
  \begin{equation*}
    \left(\ZZ^G_{\text{dim}}\left(\mu\right)\right)^2 = 
    \sum_{m_A,m_B}w^G_{\text{dim}}(m_A)w^G_{\text{dim}}(m_B).
  \end{equation*}
  
  Since the roads weights are all equal to $1$, the weight of a couple
  $(m_A,m_B)$ can be seen as a product of ``face weight'' of the
  following form; the notation $\text{d.dim}$ is a shorthand for
  ``double dimers'':
  \begin{equation}
    w_{\text{d.dim}}^f(m_A,m_B) := 
    \prod_{e\in m_A \ \mathrm{in \ city \ } f}\mu_e
    \prod_{e\in m_B \ \mathrm{in \ city \ } f}\mu_e.
  \end{equation}
  Thus,
  \begin{equation*}
    w^G_{\text{dim}}(m_A)w^G_{\text{dim}}(m_B) = 
    \prod_{f\in F} w_{\text{d.dim}}^f(m_A,m_B).
  \end{equation*}
  
  In Figure~\ref{fig:tab3}, in the first column, we have listed all
  local configurations of $(m_A,m_B)$ at a face (up to symmetries).

  \begin{figure}[!h] \centering \input{tikz/dimweightgen2.tex}
    \begin{tabular}{|c|c|c||c|c|c} \cline{1-5}
     \begin{tabular}[c]{@{}c@{}}\textbf{Dimers}\\ $m_A,
        \textcolor{green}{m_B}$ 
     \end{tabular}
   & \begin{tabular}[c]{@{}c@{}}\textbf{Fused}\\\textbf{dimers}\\$\overline{(m_A,
       \textcolor{green}{m_B})}$
     \end{tabular}
   & \begin{tabular}[c]{@{}c@{}}\textbf{Face weight}\\
       $\sum
       w_{\text{d.dim}}^f(m_A,\textcolor{green}{m_B})$ \\
       $=w_{\text{f.d.dim}}^f\left(\overline{(m_A,\textcolor{green}{m_B})}\right)$
     \end{tabular}
   & \begin{tabular}[c]{@{}c@{}}\textbf{Marked loops}\\
       $\sigma^{\text{m}}$
     \end{tabular}
   & \begin{tabular}[c]{@{}c@{}}\textbf{Face weight}\\
       $w^f_{i_{\sigma^{\text{m}}}(f)}$
     \end{tabular} & \\
      
      \cline{1-5}
      \input{tikz/dd3aa.tex} & \input{tikz/dd3aa.tex} & $a^2$ &
      \input{tikz/ll3aa.tex} & $\lambda a^2$ & * \\

      \cline{1-5}
      \input{tikz/dd3ac.tex} & \input{tikz/dd3ac.tex} & $a^2$ &
       \input{tikz/ll3ac.tex} & $\lambda a^2$ & * \\

      \cline{1-5}
      \input{tikz/dd1.tex} & \input{tikz/dd1.tex} & $1$ &
      \input{tikz/l1a.tex} \input{tikz/l1b.tex} &
      \begin{tabular}[c]{@{}c@{}}
        $\lambda a^2 + \lambda b^2$ \\ $= \lambda$
      \end{tabular} & \\

      \cline{1-5}
      \begin{tabular}[c]{@{}c@{}}
        \input{tikz/dd7a.tex} \input{tikz/dd7b.tex} \\
        \input{tikz/dd7c.tex} \input{tikz/dd7d.tex}
      \end{tabular} &
      \input{tikz/fd7.tex} &
      \begin{tabular}[c]{@{}c@{}}
        $a^4 + b^4 + 2 a^2 b^2$\\
        $=1$
      \end{tabular} &
      \input{tikz/l7a.tex} \input{tikz/l7b.tex} &
      \begin{tabular}[c]{@{}c@{}}
        $\lambda a^2 + \lambda b^2$ \\
        $= \lambda$
      \end{tabular} & \\

      \cline{1-5}
      \input{tikz/dd4.tex} &\input{tikz/dd4.tex} & $b^2$ &
      \input{tikz/l4.tex} & $\lambda b^2$ & *\\

      \cline{1-5}
      \input{tikz/dd2a.tex} & \input{tikz/dd2a.tex} & $b$ &
      \input{tikz/ll2a.tex} & $\lambda b$ & *\\

      \cline{1-5}
      \input{tikz/dd6a1.tex} \input{tikz/dd6a2.tex} & \input{tikz/fd6.tex} &
      \begin{tabular}[c]{@{}c@{}}
        $a^2b+b^3$\\
        $=b$
      \end{tabular} &
      \input{tikz/ll6a.tex} & $\lambda b$ & *\\

      \cline{1-5}
      \input{tikz/dd5.tex} & \input{tikz/dd5.tex} & $ab$ &
      \input{tikz/ll5.tex} & $\lambda a b$ & *\\
      \cline{1-5}
    \end{tabular}
    \caption{Local configurations at a face for double dimers, fused
      dimers and marked loops. For each row marked with a *, similar rows
      could be obtained by applying a symmetry or by switching the role of
      dimers $A$ and $B$.}
    \label{fig:tab3}
  \end{figure}

  \subsubsection*{Step 2: fused double-dimers} In
  Figure~\ref{fig:tab3}, the different rows correspond to the possible
  occupations of roads by double dimers, and how these roads are
  connected inside the city (up to symmetries). We can group together
  the local configurations belonging to the same row, to get a
  \textit{fused} double dimer configuration. More precisely, a fused
  dimer configuration $\overline{(m_A,m_B)}$ is the equivalence class of
  the couples of dimer configurations $(m_A,m_B)$ having identical
  roads, and having the same connections of single-dimered roads.
  
  We can define face weights for fused dimers, simply by summing the
  double dimers' face weights; they are given in the third column of
  Figure~\ref{fig:tab3}. We denote these local weights by
  $w_{\text{f.d.dim}}^f\left(\overline{(m_A,m_B)}\right)$, and the
  weight of a fused double dimer configuration is simply
  \begin{equation*}
    w^G_{\text{f.d.dim}}\left(\overline{(m_A,m_B)}\right) = 
    \prod_{f\in F} w_{\text{f.d.dim}}^f\left(\overline{(m_A,m_B)}\right).
  \end{equation*}
  
  We thus get a weight-preserving (many-to-one) mapping between a
  couple of dimer models and a fused dimer model. We represent such a
  mapping by a diagram:
  \begin{center} Dimers$_A \times $ Dimers$_B$ $\longrightarrow$ Fused
    dimers.\end{center}

  \subsubsection*{Step 3: marked loops}
  Let us now focus on the $C^{(1)}_2$ loops part. We define a
  \textit{marked} $C^{(1)}_2$ loop configuration as a $C^{(1)}_2$ loop
  configuration where every red edge is marked with an index $0$ or
  $2$, and every blue edge with an index $A$ or $B$, so that the index
  on a path stays constant except when a different-colored path is
  crossed. When a different-colored path is crossed, the index has to
  change.
  
  To make this definition well defined, we need the following Lemma:
  
  \begin{Le}
    \label{lemma:even}
    When $G$ is a connected bipartite quadrangulation on a connected
    orientable surface, for every $C^{(1)}_2$ loop configuration
    $\sigma$ on $G$, every blue loop of $\sigma$ crosses red paths an
    even number of times (and similarly for red loops).
  \end{Le}
  
  \begin{proof}
    Let us color the vertices of $G$ in black and
    white. Let us consider a blue loop of $\sigma$, and let us chose an
    orientation for that loop. When the loop turns right or left inside a
    face, the color of the vertex on its left does not change. On the
    other hand, when the loop goes straight into a face (\textit{i.e.}
    when it crosses a red loop), this color changes. Since the surface is
    orientable, when we follow the loop from start to end, we end up with
    the initial vertex on its left. Therefore it has crossed red paths an
    even number of times.
  \end{proof}
  
  For a marked loop configuration $\sigma^{\text{m}}$, let $\sigma$ be
  its unmarked version. The weight of $\sigma^{\text{m}}$ is defined
  just as in \eqref{eq:weight2}, except there is no factor for the
  number of loops; we use the notation $\text{m.loop}$ for marked loops:
  
  \begin{equation}
    \label{eq:weightmarked}
    w^G_{\text{m.loop}}(\sigma^{\text{m}}) =
    \prod_{(f,i_k)\in \sigma}w^f_i.
  \end{equation}
  
  Since there are $2^{N_{\sigma}}$ ways to mark a $C^{(1)}_2$ loop
  configuration $\sigma$, the sum of all possible marked loop weights is
  equal to the weight of $\sigma$. Therefore there is a
  weight-preserving ($2^{N_{\sigma}}$-to-one) mapping:
  \begin{center}
    Marked $C^{(1)}_2$ loops $\longrightarrow$ $C^{(1)}_2$ loops.
  \end{center}

  \subsubsection*{Step 4: marked loops to fused dimers}
  Let us describe a mapping from marked loops to fused double
  dimers. Given a marked loop configuration $\sigma^{\text{m}}$, let
  us first put dimers on roads: put one $A$ dimer on blue roads marked
  with an $A$, one $B$ dimer on blue roads marked with a $B$, two
  dimers on red roads marked with a $2$, and no dimer on red roads
  marked with a $0$. Then, let us chose the dimers in cities: when
  four blue edges come together at a face, chose the city dimers
  according to the loops connections (see the first row of
  Figure~\ref{fig:tab3}). Otherwise, there is only one possible fused
  dimer configuration with roads constructed as before. We thus
  construct a fused double dimer configuration
  $\overline{(m_A,m_B)}$. All the cases (up to symmetries) are listed
  in Figure~\ref{fig:tab3}.

  This transformation is several-to-one, indeed, loop configurations
  having four red loops marked with a $2$ (resp. a $0$) are mapped onto
  the same fused dimer configuration. However, up to a global
  multiplicative factor $\lambda$, this transformation is
  weight-preserving because, after summing over marked loops having the
  same image, the local weights of both models differ by a same factor
  $\lambda$.

  Thus there is a weight-preserving (many-to-one) mapping:
  \begin{center}
    marked $C^{(1)}_2$ loops $\longrightarrow$ Fused dimers.
  \end{center}

  \subsubsection*{Step 5: conclusion}
  To sum things up, there is a series of (many-to-one)
  weight-preserving mappings:
  \begin{center}
    Dimers$_A \times $ Dimers$_B$ $\longrightarrow$ Fused dimers
    $\longleftarrow$ Marked $C^{(1)}_2$ loops $\longrightarrow$
    $C^{(1)}_2$ loops.
  \end{center}

  This implies the equality of partition functions, which concludes
  the proof of Theorem~\ref{theo:c12dim}.
  
\end{proof}

\medskip

From there, it is natural to ask what information between dimers and
loops is kept through these transformations. At the local level, it
seems that the best connection we can get is the following:

Given a $C^{(1)}_2$ loop configuration $\sigma$, let $\sigma_b$ be the
set of blue paths in $\sigma$. If a particular set of blue paths
$\sigma^0_b$ is fixed, we say that a double dimer configuration
$(m_A,m_B)$ has paths $\sigma^0_b$ if its set of single-dimer roads is
the set of blue edges in $\sigma^0_b$, and when the four incoming
edges at a city are blue, the dimers in the city are with connection
according to $\sigma^0_b$. Then we have:
\begin{equation}
  \label{eq:locloopdim}
  \sum_{\substack{\sigma \ C^{(1)}_2 \ \mathrm{loops \ conf.} \\
      \mathrm{s.t. \ }\sigma_b=\sigma^0_b}}
  w^G_{\text{loop}}(\sigma) = \sum_{\substack{ (m_A,m_B) \
      \mathrm{dimers \ conf.} \\
      \mathrm{with \ paths\ } \sigma^0_b}}
  w^G_{\text{dim}}(m_A)w^G_{\text{dim}}(m_B).
\end{equation}

As a result, all the observables of the $C^{(1)}_2$ loop model related
to blue loops and their connectivity correspond to some observable of
the double dimer model, which can in turn be computed by determinantal
techniques \cite{Kenyon}. On the other hand, the connectivity of red
loops seems to be lost in translation. Of course, since red and blue
loops play a symmetric role, the statistics of red loops
connectivities are the same as the blue loops' and can be computed in
the same way; what we actually mean is that the joint connectivities
of red and blue loops may not be analyzed through the dimer model.

\begin{Ex}{Probability of a dimer on a road.}

  Let $e\in E^Q$ be a road in $G^Q$. We are interested in the
  probability $p$ of that road being covered by a dimer in a single
  dimer model $m_A$, with our previous dimers weights. So
  \begin{equation*}
    p = \Pr_{\text{dim}}(e \in m_A) = 
    \frac{\sum_{m_A \mid e\in m_A} w^G_{\text{dim}}(m_A)}{\ZZ^G_{\text{dim}}}.
  \end{equation*}
  When we take two independent dimer models $m_A,m_B$, the probability
  of $e$ being covered by a single dimer is then $2p(1-p)$. Because of
  the relation with loops, this is equal to the probability of $e$ being
  covered by a blue loop. Since loops are color-symmetric, this is equal
  to $\frac12$, and we deduce that
  \begin{equation*}
    p=\frac12.
  \end{equation*}

  This is true for any dimer model on a $G^Q$ with roads having weight
  $1$ and cities having weights $a,b$ such that $a^2 + b^2 = 1$. This
  property can be proven straightaway studying this dimer model, but is
  trivial in loops.
\end{Ex}

\begin{Rk}
  \label{rk:bc} We have presented the correspondence with dimers on a
  quadrangulation without boundary, but it is possible to consider, for
  instance, finite quadrangulations of the plane with a boundary. In
  that case, there is a number of external dual edges (which we think of
  as ``half-edges'', not connected together to the external face), and
  we have to chose boundary conditions for the $C^{(1)}_2$ loop
  model. For instance, we could use free boundary conditions by imposing
  nothing on the red or blue paths that use these external edges; we
  could also impose the colors of these paths to be fixed; we could even
  specify how paths starting on external edges are connected inside the
  graph.

  The equivalence with dimers works similarly by defining appropriate
  boundary conditions for dimers, except if the specification for loops
  connections concern both blue and red connections, since we cannot
  keep track of both colors' connections in our mappings.

  One example of tractable boundary conditions, where only blue
  connections are specified, will be studied in Section~\ref{sec:taut}.
\end{Rk}

\subsection{Free energy}
\label{sec:freeen}

In this paragraph, we consider an infinite quadrangulation $G$ of the
plane $\R^2$ (so it is necessarily bipartite), that is
$\Z^2$-periodic. This means that there is a basis $(e_x,e_y)$ of
$\R^2$ such that the translations by $e_x$ and $e_y$ are
color-preserving graph isomorphisms.

We define a toroidal exhaustion of $G$ in the following way: for any
$n\in \N^*$, $G_n$ is the quotient of $G$ by the lattice $n\Z e_x +
n\Z e_y$. We note $V_n,E_n,F_n$ its set of vertices, edges and
faces. For each $n$, $G_n$ is a bipartite quadrangulation on the
torus. The graph $G_1$ is called the \textit{fundamental domain} of
our quadrangulation $G$.

We assume that every face $f$ of $G_1$ is equipped with a set of
weights $w^f_1,\dots,w^f_5$ that satisfy the free-fermion
relations~\eqref{eq:ff}. Then we can define a $C^{(1)}_2$ loop model
on $G_1$ using these weights. We extend those weights periodically to
get a similar model on $G_n$: every face $f$ of $G_n$ has a unique
representative $f_0$ in $G_1$ and inherits the weights of $f_0$.

We will use the shorthand notation $\ZZ_n$ for the partition function
$\ZZ_{\text{loop}}^{G_n} \left(W\right)$ of this loop model on $G_n$.
Our goal is to compute the \textit{free energy} of this model, which
we define without a minus sign following
\cite{CohnKenyonPropp,KenyonOkounkovSheffield}:
\begin{equation*}
  \F = \lim_{n\to \infty} \frac{1}{n^2} \log\ZZ_n.
\end{equation*}
The fact that this limit exists and its exact value will follow from
the correspondence with dimers.

\medskip

We consider the dimer model of Theorem~\ref{theo:c12dim} for the
fundamental domain $G_1$.  It can be extended periodically to get a
dimer model on $G_n$. We simply note $\mu$ this set of dimers
weights. Then by Theorem~\ref{theo:c12dim},
\begin{equation*} w
  \ZZ_n = \left(\prod_{f\in F_1}\lambda_f\right)^{n^2}
  \left(\ZZ^{G_n}_{\text{dim}}\left(\mu \right)\right)^2
\end{equation*}
so that
\begin{equation}
  \label{eq:zldimlog}
  \frac{1}{n^2} \log\ZZ_n = \sum_{f\in F_1}\log(\lambda_f) + \frac{2}{n^2}
  \log\left(\ZZ^{G_n}_{\text{dim}}\left(\mu\right)\right).
\end{equation}
The right-hand side of \eqref{eq:zldimlog} contains the free energy of
the periodic dimer model on $G^Q$ with weights $\mu$. This quantity
can be exactly computed
\cite{CohnKenyonPropp,KenyonOkounkovSheffield}. Let us recall how this
computation is made.

\medskip

The graph $G^Q_1$, is a bipartite graph on the torus. We equip it with
a Kasteleyn orientation (\textit{i.e.} an orientation of edges such
that every face has an odd number of clockwise edges; there exists
such an orientation, see for instance \cite{CimasoniReshetikhin}). We
can split its vertices between black ones $B^Q_1$ and white ones
$W^Q_1$. Then we define its Kasteleyn matrix $K_1$ as a $W^Q_1$ by
$B^Q_1$ weighted adjacency matrix, with entries
\begin{equation*}
  (K_1)_{w,b} =
  \begin{cases}
    \mu(e) & \text{when} \begin{tikzpicture}
    \draw [] (0,0) -- node [above] {$e$} (1,0); 
    \draw [->, >=stealth] (0.5,0) -- (0.52,0); 
    \node [draw=black, fill=white,thick,circle,inner sep=0pt,minimum size=5pt] (a1) at (0,0) {};
    \node [draw=black, fill=black,thick,circle,inner sep=0pt,minimum size=5pt] (a2) at (1,0) {};
    \draw (a1) node [above] {$w$};
    \draw (a2) node [above] {$b$};
\end{tikzpicture} \\
    -\mu(e) & \text{when} \begin{tikzpicture}
    \draw [] (0,0) -- node [above] {$e$} (1,0); 
    \draw [->, >=stealth] (0.5,0) -- (0.48,0); 
    \node [draw=black, fill=white,thick,circle,inner sep=0pt,minimum size=5pt] (a1) at (0,0) {};
    \node [draw=black, fill=black,thick,circle,inner sep=0pt,minimum size=5pt] (a2) at (1,0) {};
    \draw (a1) node [above] {$w$};
    \draw (a2) node [above] {$b$};
\end{tikzpicture} \\
    0 & \text{otherwise.}
  \end{cases}
\end{equation*}
Let $z$ and $w$ be two complex numbers. We construct a modified
matrix $K_1(z,w)$ in the following way.  Let $\gamma_x,\gamma_y$ be
the two oriented cycles on the torus corresponding respectively to
$e_x,e_y$. We multiply the weight of edges crossing $\gamma_x$ by $z$
when the white vertex is on the left of $\gamma_x$, and by $z^{-1}$
when the white vertex is on the right of $\gamma_x$; and similarly for
$\gamma_y$ and $w$. These weights define a Kasteleyn matrix
$K_1(z,w)$. The \textit{characteristic polynomial} of $G_1$ is then:
\begin{equation*}
  P(z,w) = \det K_1(z,w).
\end{equation*}

The free energy of the dimer model can be expressed using the
characteristic polynomial:
\begin{Theo*}[\cite{CohnKenyonPropp,KenyonOkounkovSheffield}]
  We have
  \begin{equation*}
    \lim_{n\to\infty} \frac{1}{n^2}
    \log\left(\ZZ^{G_n}_{\text{dim}}\left(\mu\right)\right) =
    \frac{1}{(2i\pi)^2}\int_{\T^2}\log |P(z,w)| \frac{dz}{z} \frac{dw}{w}.
  \end{equation*}
\end{Theo*}

Therefore, \eqref{eq:zldimlog} gives
\begin{equation}
  \F = \sum_{f\in F_1}\log(\lambda_f) -
  \frac{1}{2\pi^2}\int_{\T^2}\log |P(z,w)| \frac{dz}{z}
  \frac{dw}{w}.
\end{equation}

\begin{Ex}{Free energy for integrable weights on $\Z^2$.}
  Let us take for $G$ the quadrangulation $\Z^2$.  We equip every face
  with the positive integrable weights of \cite{WarnaarNienhuis,
    IkhlefCardy} for the fugacity $n=2$:
  \begin{equation}
    \begin{cases}
      w_1 = \sin^2 (\theta) \\
      w_2 = \cos^2 (\theta)\\
      w_3 = \sin (\theta)\\
      w_4 = \cos (\theta)\\
      w_5 = \sin (\theta) \cos (\theta)
    \end{cases}
  \end{equation}
  where $\theta \in (0,\frac{\pi}{2})$. Then for all
  $f\in F$, $\lambda_f=1,a_f=\sin(\theta),b_f=\cos(\theta)$.
  
  The corresponding weighted graph $G^Q_1$, equipped with a Kasteleyn
  orientation, is represented on the left of Figure~\ref{fig:octa}. We
  can perform an urban renewal \cite{KenyonProppWilson} on this graph to
  get the graph on the right.
  
  \begin{figure}[!h]
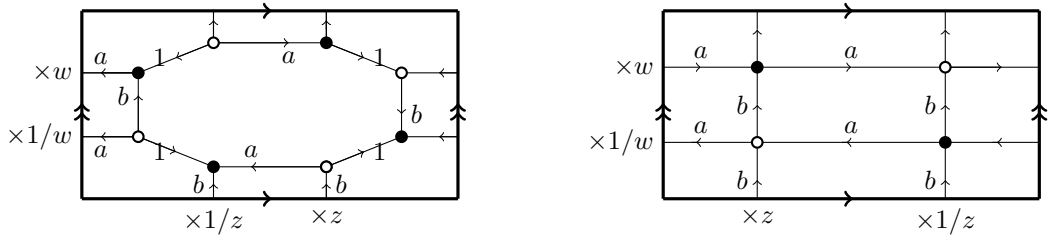

    \centering \input{tikz/octa.tex}
    \hspace*{1cm}
    \centering \input{tikz/octa2.tex}
    \caption{The dimers' fundamental domain $G^Q_1$, before (left) and
      after (right) urban renewal, equipped with a Kasteleyn
      orientation. Here $a=\sin(\theta)$, $b=\cos(\theta)$.}
    \label{fig:octa}
  \end{figure}

  By using the free energy for dominoes on $\Z^2$ \cite{Fisher} with
  horizontal weight $\sin(\theta)$ and vertical weight $\cos(\theta)$,
  we get
  \begin{equation*}
    \F = - \frac{1}{2\pi^2}\int_{\T^2}\log \left|- 2 +
      \cos^2(\theta) \left(z+\frac1z\right)+ \sin^2(\theta)
      \left(w+\frac1w\right)\right| \frac{dz}{z} \frac{dw}{w}.
  \end{equation*}
  Several other expressions can be given for this
  quantity, for instance following \cite{Kenyon02}:
  \begin{equation*}
    \F = \frac{2}{\pi} L\left(\theta\right) +
    \frac{2}{\pi} L\left(\frac{\pi}{2} - \theta\right) +
    \frac{2\theta}{\pi} \ln(\tan(\theta)) + \ln(2\cos(\theta))
  \end{equation*}
  where $L$ is the Lobachevsky function.

\end{Ex}

\section{Cube flip and Kashaev's recurrence}
\label{sec:st}

From now on, $G$ might be a planar quadrangulation in the sense of
Definition~\ref{def:quadrangulation}, or a finite planar
quadrangulation with a boundary, in the following sense:

\begin{Def}
  \label{def:quadbound} A quadrangulation with a boundary is a finite
  simple graph $G=(V,E)$, properly embedded in the plane, such that all
  \textit{internal} (bounded) faces have degree $4$.
\end{Def}

In this case, we define a $C^{(1)}_2$ loop model on $G$ by specifying
a boundary condition in any way discussed in Remark~\ref{rk:bc} --
since the correspondence with dimers won't be used, both blue and red
connections can be specified.

We start by proving a few identities on Kashaev's relation, then we
relate the star-triangle relation for loops to Kashaev's relation.

\subsection{Kashaev's recurrence}

Kashaev's relation reads:
\begin{equation}
  \label{eq:kasharel}
  \begin{split}
     g^2 g_{123}^2 & + g_{1}^2 g_{23}^2 + g_{2}^2 g_{13}^2 + g_{3}^2
     g_{12}^2\\
     & - 2 g_{2}g_{3}g_{13}g_{12} - 2 g_{1}g_{3}g_{23}g_{12} -2
     g_{1}g_{2}g_{23}g_{13} \\
     & - 2 g g_{123} (g_{1} g_{23} + g_{2} g_{13} + g_{3} g_{12}) \\
     & - 4 g g_{12} g_{23} g_{13} - 4 g_{123} g_1 g_2 g_3 = 0.
  \end{split}
\end{equation}

If seven of the $g$ variables are positive real numbers, then the
eighth one can be deduced from \eqref{eq:kasharel}, up to the choice
of the square root of a quadratic polynomial.  The choice of the
greatest root guarantees that the $g$ variables stay positive. In this
case, the recurrence can be explicitly written, see Proposition
\ref{prop:calculs}.

\begin{figure}[h] \centering \input{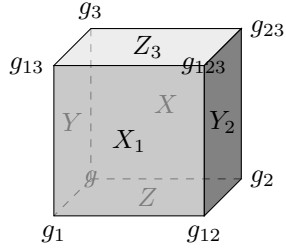}
  \caption{The $g$, $X$, $Y$ and $Z$ variables, implicitly taken at
    $x\in \Z^3$.}
  \label{fig:var}
\end{figure}

We define six other variables $X,Y,Z,X_1,Y_2,Z_3$ by
\begin{equation}
  \begin{split}
    X & =  \sqrt{g g_{23} + g_2 g_3}, \\
    Y & =  \sqrt{g g_{13} + g_1 g_3}, \\
    Z & =  \sqrt{g g_{12} + g_1 g_2}, \\
    X_1 & =  \sqrt{g_1 g_{123} + g_{12} g_{13}}, \\
    Y_2 & =  \sqrt{g_2 g_{123} + g_{12} g_{23}}, \\
    Z_3 & =  \sqrt{g_3 g_{123} + g_{13} g_{23}}.
  \end{split}
\end{equation}
All of these quantities can be nicely represented on
the vertices and faces of a cube, see Figure~\ref{fig:var}.

The following relations can be obtained by simple calculations:
\begin{Prop}[Kashaev's recurrence]
  \label{prop:calculs}
  When Kashaev's relation \eqref{eq:kasharel} is satisfied, the value
  of $g_{123}$ obtained by taking the greatest root reads
  \begin{enumerate}
  \item \label{f123} $g_{123} = \frac{2g_1 g_2 g_3 + g(g_1 g_{23} +
      g_2 g_{13} + g_3 g_{12}) +2XYZ}{g^2}.$
  \end{enumerate}
  Furthermore,
  \begin{enumerate} \setcounter{enumi}{1}
  \item \label{fx} $X_1 = \frac{g_1 X + Y Z}{g},$
  \item \label{fy} $Y_2 = \frac{g_2 Y + X Z}{g},$
  \item \label{fz} $Z_3 = \frac{g_3 Z + X Y}{g},$
  \item \label{ff} $\frac{X_1 Y_2 Z_3 + g_{12}g_{13}g_{23}}{g_{123}} =
    \frac{X Y Z + g_{1}g_{2}g_{3}}{g}.$
  \end{enumerate}
\end{Prop}
The first 4 relations are a particular case of the \textit{hexahedron
  recurrence} \cite{KenyonPemantle}.

Kashaev's relation \eqref{eq:kasharel} has the same symmetries as the
cube. In particular, by considering the central symmetry relative to
the center of the cube, we get that the previous transformation is
self-dual, in the following sense:
\begin{Prop}\label{prop:dual}
  Let $g_{123}$ be defined in terms of
  $g,g_1, g_2, g_3, g_{12}, g_{13}, g_{23}$ as in Item \ref{f123} of
  Proposition~\ref{prop:calculs}. Then the transformation:
  \begin{equation*}
    \begin{pmatrix}
      g\\ g_1\\ g_2 \\ g_3 \\ g_{23} \\ g_{13} \\ g_{12}
    \end{pmatrix}
    \mapsto
    \begin{pmatrix}
      g_{123} \\ g_{23} \\ g_{13} \\ g_{12}\\ g_1 \\ g_2 \\ g_3
    \end{pmatrix}.
  \end{equation*}
  is an involution.
\end{Prop}
Also note that this involution exchanges $X$ and $X_1$,
$Y$ and $Y_2$, $Z$ and $Z_3$.

\subsection{Parametrization of a free-fermionic $C^{(1)}_2$ loop
  model}

Let us suppose that there is a set of positive real values on the
vertices $V$ of $G$, denoted $g=(g_x)_{x\in V}$, such that the local
loop weights take the following form at a face
$f =\raisebox{-0.4cm}{\input{tikz/face3.tex}}$:

\begin{equation}
  \label{eq:paramg}
  \begin{cases}
    \raisebox{-0.15cm}{\input{tikz/a.tex}} \hspace*{0.5cm} w^f_1 = g_xg_y
    \\
    \raisebox{-0.15cm}{\input{tikz/b.tex}} \hspace*{0.5cm} w^f_2 = g_u g_v
    \\
    \raisebox{-0.15cm}{\input{tikz/c.tex}} \hspace*{0.5cm}
    w^f_3 = \sqrt{g_xg_y} \sqrt{g_xg_y + g_ug_v} \\
    \raisebox{-0.15cm}{\input{tikz/d.tex}} \hspace*{0.5cm} w^f_4 =
    \sqrt{g_vg_u} \sqrt{g_xg_y + g_ug_v}\\
    \raisebox{-0.15cm}{\input{tikz/e.tex}} \hspace*{0.5cm} w^f_5 =
    \sqrt{g_xg_yg_ug_v}.
  \end{cases}
\end{equation}

Notice that these weights satisfy the free-fermionic relations
\eqref{eq:ff}. In fact, on a class of planar quadrangulations that
includes finite lozenge graphs (embedded graphs whose internal faces
are non-degenerate rhombi with same edge length), every free-fermionic
$C^{(1)}_2$ loop model can be parametrized in this way. This is proved
in Appendix~\ref{app:fullc12}.

Notice also that the bipartite coloring of $G$ is no longer important
to define the weights. We will no longer show this coloring.

We can define a marked model with these weights using the weights of
equation~\eqref{eq:weightmarked}. We will show that the Yang-Baxter
equations for this marked loop model are equivalent to $g$ satisfying
Kashaev's relation. This implies a similar statement for non-marked
loops as well as for the double dimer model. However, since the
indices on loops don't affect the weight, the proof for marked loops
will not be any more difficult that what it would be for non-marked
loops.

\medskip

We now denote by $\ZZ^G_{\text{m.loop}}(g) $ the partition function of
marked loops on $G$:
\begin{equation*}
  \ZZ^G_{\text{m.loop}}(g) =
  \sum_{\sigma} w^G_{\text{m.loop}}(\sigma)
\end{equation*}
where $\sigma$ runs over all marked $C^{(1)}_2$ loop configurations on
$G$, and $w^G_{\text{m.loop}}(\sigma)$ is the same as
\eqref{eq:weightmarked} with the local weights defined by
\eqref{eq:paramg}.

\subsection{Cube flip}

\begin{figure}[h] \centering
  \input{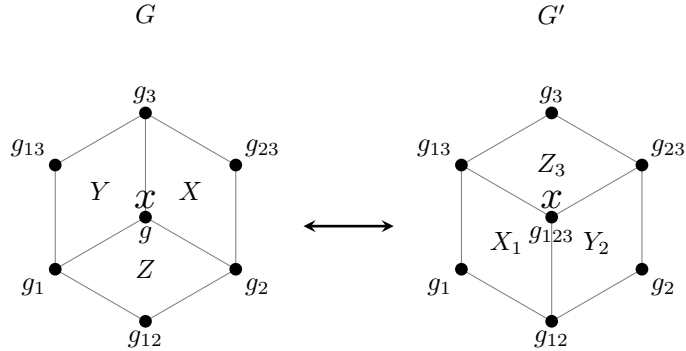}
  \caption{Cube flip at $x$.}
  \label{fig:st}
\end{figure}

Let us suppose that $G$ contains a vertex $x\in V$ of degree $3$, such
that the graph around $x$ looks like the left-hand side of
Figure~\ref{fig:st}. We can perform a ``cube flip'' at $x$ by changing
the edges around this vertex. It gives a new graph $G'$, on the
right-hand side of Figure~\ref{fig:st}. This is a form of
star-triangle (or $Y-\Delta$) move.

When such a move is performed, we allow $g_x$ to change into a new
value $g'_x$; let $g'$ be the set of variables equal to $g$ everywhere
except at $x$ where it is $g'_x$.

By labeling the $g$ variables around $x$ as in Figure~\ref{fig:st}, we
say that $g$ \textit{satisfies Kashaev's recurrence} at $x$ when
$g,g_1,\dots,g_{123}$ satisfy \eqref{eq:krec}. Note that performing
cube flips at $x$ twice with $g$ satisfying Kashaev's recurrence would
bring back the original graph and constants, because of
Proposition~\ref{prop:dual}.

\begin{Theo}
  \label{theo:cuberem}
  When $g$ satisfies Kashaev's recurrence at $x$,
  \begin{equation*}
    \frac{1}{g_x^2}\ZZ^G_{\text{m.loop}}(g) = 
    \frac{1}{(g'_x)^2}\ZZ^{G'}_{\text{m.loop}}(g').
  \end{equation*}
\end{Theo}

Another way to phrase this is to say that the Yang-Baxter equations
for the $C^{(1)}_2$ loop model, as defined in \cite{WarnaarNienhuis},
taken in our parametrization, become equivalent to Kashaev's
recurrence.

\begin{proof}
  There are six dual edges incoming in the region at $x$. We call
  \textit{connection pattern} any way to color, label and connect these
  six incoming edges. For example, in the first row of
  Figure~\ref{fig:mapping2}, the connection pattern is the following:
  all edges are all colored in blue; the west and northwest edges are
  labeled $k$ and connected; the east and northeast edges are labeled
  $l$ and connected; the southeast and southwest edges are labeled $m$
  and connected.

  We want to construct a coupling between marked $C^{(1)}_2$ loops on
  $G$ and on $G'$ so that they agree everywhere except inside the
  changed region around $x$, and such that the connection pattern
  doesn't change. Thus we want to group the marked loop configurations
  on $G$ and on $G'$ according to their connection pattern. All possible
  cases are listed in Figure~\ref{fig:mapping2}; indices on loops are
  arbitrary and we used the notation $\hat{k}$ for the index different
  from $k$ (so $\hat{A}=B$, $\hat{B}=A$, $\hat{0}=2$, $\hat{2}=0$). We
  omitted a few cases: first, the ones that can be derived from
  represented ones by a simple rotation, symmetry or color swap;
  secondly, for any row $i\in \{1,\dots, 7\}$, there is a dual row
  obtained by taking the right-hand side configurations of $i$, turning
  them upside-down and drawing them on the left - and similarly for the
  other side (notice that rows $2,6,7$ are self-dual).

  \begin{figure}[!h]
    \centering
    \input{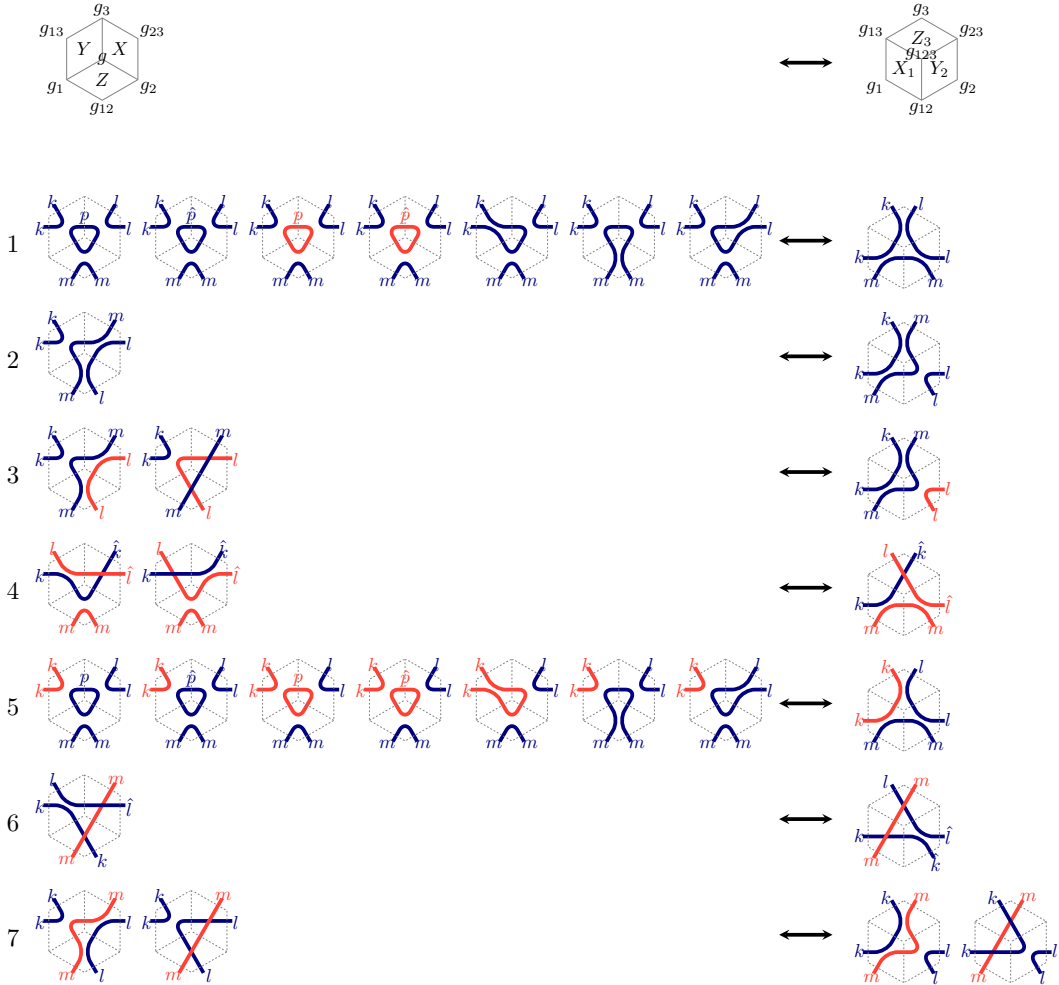}
    \caption{Correspondence between marked loops on $G$ (left) and
      $G'$ (right) having the same connection pattern.}
    \label{fig:mapping2}
  \end{figure}

  Note that connection patterns are the same for $G$ and $G'$. Let $m$
  be the total number of connection patterns, and let $\{1,\dots, m\}$
  be a set of indices representing them. For any $j\in \{1,\dots, m\}$,
  let $\Sigma_j$ be the set of marked $C^{(1)}_2$ loop configurations on
  $G$ having that connection pattern. Then $\Sigma_1, \dots, \Sigma_m$
  are a partition of the set of all marked configurations on $G$. Thus
  in terms of partition functions we have
  \begin{equation*}
    \ZZ^G_{\text{m.loop}}(g) = \sum_{j=1}^{m} \ZZ_j
  \end{equation*}
  where $\ZZ_j = \sum_{\sigma \in \Sigma_j}
  w^G_{\text{m.loop}}(\sigma)$.  Similarly,
  \begin{equation*}
    \ZZ^{G'}_{\text{m.loop}}(g) = \sum_{j=1}^{m} \ZZ'_j
  \end{equation*}
  where $\ZZ'_j = \sum_{\sigma \in \Sigma'_j}
  w_{\text{m.loop}}^{G'}(\sigma)$, and $\Sigma'_j$ being the set of
  marked $C^{(1)}_2$ loop configurations on $G'$ having connection
  pattern $j$.

  If $\sigma \in \Sigma_j$, then its total weight can be written as

  \begin{equation*}
    w^G_{\text{m.loop}}(\sigma) = a(\sigma) w_{\text{loc.}} (\sigma)
  \end{equation*}
  where $w_{\text{loc.}}(\sigma) = w^X_{i_1}w^Y_{i_2}w^Z_{i_2}$, with
  $(X,i_1),(Y,i_2),(Z,i_3) \in \sigma$; this is the local weight coming
  from the faces around $x$, and $a(\sigma)$ doesn't depend on the local
  configuration around $x$. Actually, if two local configurations have
  the same connection pattern, the possible ways to extend them to
  construct a loop configuration are the same, so the possible values of
  $a(\sigma)$ for these configurations are the same. This shows that
  $\ZZ_j$ can be factored into:
  \begin{equation}
    \label{eq:zj}
    \ZZ_j = A_j \left( \sum_{\sigma \in L_j}
      w_{\text{loc.} }(\sigma) \right),
  \end{equation}
  where $A_j$ is the sum of weight of all possible
  $a(\sigma)$ for any $\sigma \in \Sigma_j$, and $L_j$ is the set of all
  local configurations around $x$ in $G$ that have connection pattern
  $j$.

  Similarly, $\ZZ'_j$ takes the form
  \begin{equation}
    \label{eq:zpj}
    \ZZ'_j = A_j \left(\sum_{\sigma \in L'_j}
      w_{\text{loc.} }(\sigma)\right).
  \end{equation}

  Notice that the part $A_j$ is the same for $\ZZ_j$ and $\ZZ'_j$:
  indeed, $G$ and $G'$ are the same outside of the region around $x$;
  given a local configuration in that region, the list of possible
  extensions into a global loop configuration only depends on its
  connection pattern.

  As a result, the following lemma is sufficient to conclude the
  proof.

  \begin{Le}
    \label{lemma:map}
    For any row $i \in \{1,\dots,7\}$ in Figure~\ref{fig:mapping2},
    let $L_i$ (resp. $L'_i$) be the list of local configurations in
    row $i$ on $G$ (resp. $G'$). If $g$ satisfies Kashaev's
    recurrence, then
    \begin{equation*}
      \frac{1}{g_x^2}\sum_{\sigma \in L_i} w_{\text{loc.} }(\sigma) =
      \frac{1}{(g'_x)^2}\sum_{\sigma \in L'_i} w_{\text{loc.} }(\sigma).
    \end{equation*}
  \end{Le}

  \begin{proof}
    For $i=1$, $L_1$ contains 7 configurations, and the
    sum of their local weights is, in the same order as in
    Figure~\ref{fig:mapping2} (using the notations $g,g_1,\dots, g_{123}$
    like in Figure~\ref{fig:st}),
    \begin{equation}
      \label{eq:ra}
      \begin{split}
        \frac{1}{g_x^2} \sum_{\sigma \in L_1} w_{\text{loc.} }(\sigma)
        & = \frac{2g_1^2 g_2^2 g_3^2}{ g^2} + \frac{2g_1 g_2 g_3 X Y
          Z}{ g^2} + \frac{g_1 g_2^2 g_3 g_{13}}{ g} + \frac{g_1 g_2
          g_3^2 g_{12}}{ g} + \frac{g_1^2 g_2 g_3 g_{23}}{g} \\
        & = g_1 g_2 g_3 \left( \frac{2g_1 g_2 g_3 + g(g_1 g_{23} + g_2
            g_{13} + g_3 g_{12})+2XYZ}{g^2} \right).
      \end{split}
    \end{equation}

    There is only one configuration in $L'_1$ which has local weight:
    \begin{equation}
      \label{eq:la}
      \frac{1}{(g'_x)^2}\sum_{\sigma \in
        L'_1} w_{\text{loc.} }(\sigma) = g_1 g_2 g_3 g_{123}.
    \end{equation}
    The equality of \eqref{eq:la} and \eqref{eq:ra} is given by Item
    \ref{f123} in Proposition~\ref{prop:calculs}.

    The other cases are similar, using the various relations of
    Proposition~\ref{prop:calculs}. This is done in
    Appendix~\ref{app:map}.
  \end{proof}

  Because of the self-duality property (Proposition~\ref{prop:dual})
  and the symmetries of the model, Lemma~\ref{lemma:map} is also true
  for any represented or non-represented row. Since connection patterns
  of configurations on $G$ and $G'$ are the same, the boundary
  conditions (if any) are preserved by this coupling. This concludes the
  proof of Theorem~\ref{theo:cuberem}
\end{proof}

\begin{Rk}
  Based on Theorem~\ref{theo:cuberem}, it is natural to define a
  renormalized partition function
  \begin{equation*}
    \YY^G_{\text{m.loop}}(g) =
    \left( \prod_{x\in V} \frac{1}{g_x^2}\right)
    \ZZ^G_{\text{m.loop}}(g).
  \end{equation*}
  
  We can also go back to unmarked loops, weighted by a factor
  $2^{N_{\sigma}}$ \eqref{eq:weight2}, so that
  \begin{equation}
    \label{eq:yzinv}
    \YY^G_{\text{loop}}(g) =
    \YY^{G'}_{\text{loop}}(g').
  \end{equation}
  This quantity will appear again as the combinatorial object
  representing the solution of Kashaev's recurrence on a stepped
  surface.
\end{Rk}

\begin{Rk}
  In the language of statistical physics, Theorem~\ref{theo:cuberem}
  is a case of Z-invariance.  Since marked loops generalize unmarked
  loops as well as dimers, $Z$-invariance for marked loops under
  Kashaev's recurrence implies $Z$-invariance for unmarked loops, and
  (when the boundary conditions, if any, don't involve red
  connections) for the double dimer model of Theorem~\ref{theo:c12dim}
  for parametrized weights \eqref{eq:paramg}.

  This particular dimer model is represented in
  Figure~\ref{fig:dimkp}. After performing gauge transformations (by
  multiplying weights by $\frac{1}{\sqrt{g_xg_u}}$ around any vertex
  of $G^Q$ that is closest to the edge $\{xu\}$ of $G$), we get the
  weights on the right of Figure~\ref{fig:dimkp}. This is the dimer
  model of \cite{KenyonPemantle} in the particular case of Kashaev's
  relation; we get an alternative proof of its Z-invariance.

  \begin{figure}[!h]
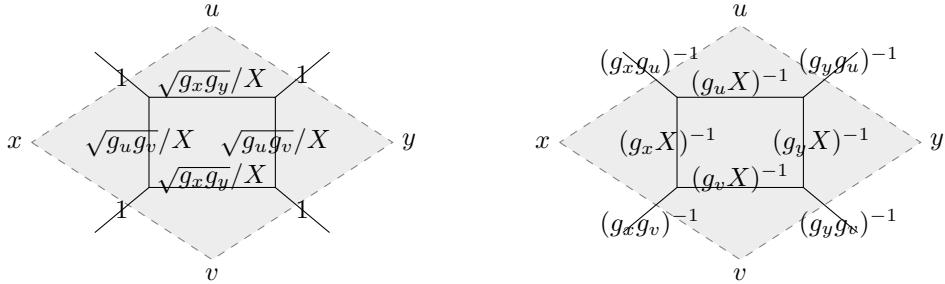

    \centering
    \input{tikz/dimweight2.tex}
    \hspace*{1cm}
    \input{tikz/dimweight.tex}
    \caption{The dimer model's weights before (left) and after (right)
      gauge transformations. Here we have set $X=\sqrt{g_x g_y + g_u g_v}$.}
    \label{fig:dimkp}
  \end{figure}

  Interestingly enough, if one starts from Kashaev's parametrized
  Ising model \eqref{eq:kashis}, and applies the procedure transforming
  Ising models into dimers on $G^Q$ \cite{Dubedat,BoutillierDeTiliere},
  one gets the same dimer model as in Figure~\ref{fig:dimkp}. However,
  the transformation between Ising model and dimers is not a direct
  mapping, so it is not completely clear in that setting why the
  Yang-Baxter equations for the Ising model translate into the same
  equations for dimers.
\end{Rk}

\section{Taut configurations on stepped surfaces}
\label{sec:taut}

We now turn to the study of \textit{taut} configurations, which will
be the appropriate objects counted by the solution of Kashaev's
recurrence for arbitrary initial conditions.

\subsection{Stepped surfaces}

We denote by $(e_1,e_2,e_3)$ the canonical basis of $\R^3$.  For
$x=(i,j,k)\in \Z^3$, let $C_{x}\subset \R^3$ be the unit cube
$[i,i+1]\times [j,j+1] \times [k,k+1]$.

We call \emph{stepped solid} a union of such unit cubes. A stepped
solid $U$ is said to be \emph{monotone} if, for every $C_{(i,j,k)}
\subset U$, and for every $i'\leq i, j' \leq j, k' \leq k$,
$C_{(i',j',k')} \subset U$.

In this section, we always assume that $U$ is a monotone stepped
solid. In that case, the topological boundary $\partial U$ is a union
of squares of the form $(x,x+e_i,x+e_i+e_j,x+e_j)$ where $x\in \Z^3$
and $1\leq i < j \leq 3$; it is called a \emph{stepped surface}. This
boundary naturally corresponds to an infinite planar quadrangulation
$G(U)$, formally defined by the following sets of vertices and edges:
\begin{equation*}
   V(U) = \{x\in U\cap \Z^3 \mid x+e_1+e_2+e_3 \notin U\},
\end{equation*}

\begin{equation*}
  E(U) = \left\{\{x,y\} \mid x,y \in V, \  x-y \in 
  \{\pm e_1, \pm e_2, \pm e_3\} \right\}.
\end{equation*}

We will also denote $F(U)$ the set of faces of $G(U)$.

\begin{figure}[h]
  \centering
  \input{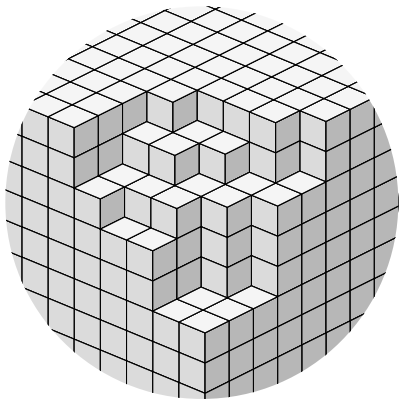}
  \caption{A portion of a stepped solid $U$.}
  \label{fig:exstep}
\end{figure}

\subsection{Boundary conditions}
To define a loop model on $G(U)$ we need some extra boundary
conditions. Consider the negative corner
$U_0 = \R_- \times \R_- \times \R_-$. We require that $U \subset U_0$
and that $U_0 \setminus U$ be composed of finitely many cubes. In this
case, we say that $U$ is \emph{regular} (see Figure~\ref{fig:exstep}
for example). The infinite graph $G(U)$ is then equal to $G(U_0)$
outside of a sufficiently big ball for the euclidean distance centered
at the origin. Let $B(O,R_U)$ be such a ball.

If $U$ is a regular stepped solid, and $(g_x)_{x\in V(U)}$ is a
collection of variables on the vertices of $G(U)$, then Kashaev's
recurrence \eqref{eq:krec} is enough to recursively define $g$ on
$\Z_-^3 \setminus U$. Then, the value at the origin $g_{(0,0,0)}$ is
called the \textit{solution of Kashaev's recurrence} with initial
conditions $(g_x)_{x \in V(U)}$.

\medskip

\begin{figure}[h]
  \centering
  \input{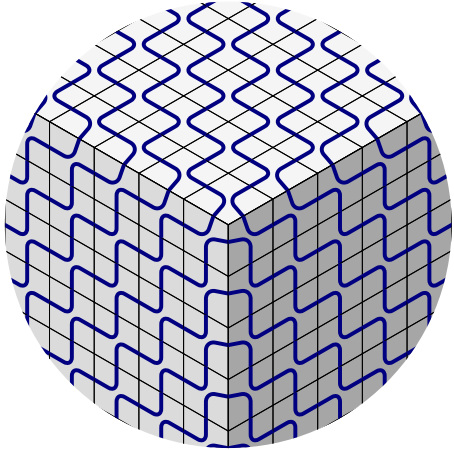}
  \caption{The initial stepped solid $U_0$ and initial configuration
    $\sigma_0$.}
  \label{fig:zero}
\end{figure}

On $G(U_0)$ let $\sigma_0$ be the configuration given in Figure
\ref{fig:zero}. Following the terminology of \cite{KenyonPemantle}, we
say that a $C^{(1)}_2 $ loop configuration $\sigma$ on $G(U)$ is
\emph{taut} when it has the same connectivity as $\sigma_0$ on a
neighborhood of infinity -- meaning that outside of a ball of radius
$R_{\sigma} \geq R_U$, $\sigma$ has to be equal to $\sigma_0$, and
$\sigma$ has to contain paths connecting the edges of $U\cap
B(O,R_{\sigma})$ that are connected in $\sigma_0$. This is the case in
Figure~\ref{fig:ex} for example.

Let $\Sigma(U)$ be the set of taut $C^{(1)}_2$ loop configuration on
$G(U)$.

The following lemmas are direct adaptations of Carroll and Speyer's
arguments on cube groves \cite{CarrollSpeyer}.

\begin{Le}
  \label{lemma:ci}
  $\Sigma(U_0)$ only contains $\sigma_0$.
\end{Le}

\begin{proof}
  Let $\sigma$ be a taut configuration on $G(U_0)$. Suppose that it is
  different from $\sigma_0$. Without loss of generality, we can assume
  that there is a face $f$ of the form $\{x,x+e_1,x+e_1+e_2,x+e_2\}$
  that differs from $\sigma_0$. Take such a face with $x=(i,j,0)$ and
  $i+j$ minimal; this is possible because $\sigma$ is equal to
  $\sigma_0$ for all faces far enough from the origin. Then the faces
  $f-e_1$ and $f-e_2$ have to be the same for $\sigma$ and $\sigma_0$,
  which implies that the dual edges of $\{x,x+e_1\}$ and $\{x,x+e_2\}$
  are blue. Now since $\sigma$ differs from $\sigma_0$ on $f$, the
  only possible local configurations connect these two blue edges
  inside of $f$. This implies that the connectivity of $\sigma$
  differs from the one in $\sigma_0$, which contradicts the definition
  of a taut configuration.
\end{proof}

Note that a configuration $\sigma$ contains two types of paths:
infinite simple paths, and finite closed simple loops. In the rest of
the paper, the latter will simply be called ``loops''.

\begin{Le}
  \label{lemma:fini}
  For a regular stepped solid $U$, $\Sigma(U)$ is finite. Moreover, a
  taut $C^{(1)}_2$ loop configuration on $U$ has a finite number of
  loops.
\end{Le}

\begin{proof}
  Let $\sigma \in \Sigma(U)$. We know that $\sigma$ has the same
  connectivity as $\sigma_0$ and is equal to $\sigma_0$ outside of
  $B(O,R_{\sigma})$. By the same argument described in the proof of
  Lemma~\ref{lemma:ci}, $\sigma$ actually has to be equal to
  $\sigma_0$ outside of $B(O,R_U)$. For a fixed $U$, there is only a
  finite number of such $\sigma$, and their loops have to be in
  $B(O,R_U)$ so there is a finite number of them.
\end{proof}

\subsection{Weights}
\label{sweights}

Let $U$ be a regular stepped solid, and let $(g_x)_{x \in V(U)}$ be a
collection of variables on the vertices of $G(U)$ that can be thought
of as positive real numbers.  For a taut $C^{(1)}_2$ loop
configuration $\sigma \in \Sigma(U)$, we still denote $N_{\sigma}$ its
number of loops. Let us consider the weight:
\begin{equation}
  \label{eq:weight}
  w_{\text{taut}}^U(\sigma) =
  2^{N_{\sigma}} \left( \prod_{(f,i_k)\in \sigma}w^f_{i} \right)
  \left( \prod_{x\in V(U)} \frac{1}{g_x^2}\right),
\end{equation}
where the local weights $w^f_i$ are defined using parametrization
\eqref{eq:paramg}. This expression makes sense because $N_{\sigma}$ is
finite (Lemma~\ref{lemma:fini}) and $\prod w^f_{i}$ formally makes
$g_x$ appear with exponent $2$ for every $x$ in the ``flat'' region
where $\sigma = \sigma_0$, so the two products cancel out for all but
a finite number of terms. For example,
$w_{\text{taut}}^{U_0}(\sigma_0)=g_{(0,0,0)}$.

Since the number of taut configurations is finite
(Lemma~\ref{lemma:fini}), we can define the renormalized partition
function:
\begin{equation*}
  \YY^U_{\text{taut}}(g) = \sum_{\sigma \in \Sigma(U)}
  w_{\text{taut}}^U(\sigma).
\end{equation*}

Formally, $\YY^U_{\text{taut}}(g)$ is simply a Laurent polynomial in
the $g$ variables and in the $\sqrt{g_x g_y + g_u g_v}$ variables, for
any face bordered by $x,u,y,v$; we call these \textit{face variables}.

\begin{Theo}
  \label{theo:princ2}
  Let $U$ be a regular stepped solid, and $(g_x)_{x\in V(U)}$ be a
  collection of variables on the vertices of $G(U)$. Then the solution
  of Kashaev's recurrence at the origin with initial condition
  $(g_x)_{x \in V(U)}$ is
  \begin{equation*}
    g_{(0,0,0)} = \YY^U_{\text{taut}}(g).
  \end{equation*}
\end{Theo}

\begin{proof}
  Let $g$ be the solution of Kashaev's recurrence with initial
  conditions on $V(U)$ (so that $g$ is also defined on
  $\Z_-^3\setminus U$).

  The results of Section~\ref{sec:st} imply that the renormalized
  partition function $\YY^U_{\text{taut}}(g)$ doesn't change when a
  cube is added to, or removed from, $U$ in such a way that it remains
  a stepped solid - as long as we keep using the $g$
  variables. Indeed, the boundary conditions are unchanged when a cube
  flip is performed so all the computations stay the same for taut
  configurations.

  By repeatedly removing cubes starting from $U_0$ to get to $U$, this
  implies
  \begin{equation*}
    \YY^U_{\text{taut}}(g) = \YY^{U_0}_{\text{taut}}(g).
  \end{equation*}

  But $\Sigma(U_0)$ only contains $\sigma_0$ (Lemma~\ref{lemma:ci}),
  and $w_{\text{taut}}^{U_0}(\sigma_0)=g_{(0,0,0)}$, so that
  $\YY^{U_0}_{\text{taut}}(g) = g_{(0,0,0)}$.
\end{proof}

\subsection{Algebraic consequences}

We have seen that the solution of Kashaev's recurrence is equal to the
partition function $\YY^U_{\text{taut}}(g)$.  However, several taut
$C^{(1)}_2$ loop configurations on $G(U)$ might correspond to the same
monomial in $\YY^U_{\text{taut}}(g)$. The following Theorem states
that it is not the case, and gives consequences on the exponents and
coefficients appearing in the Laurent polynomial. The first and third
points were already obtained in \cite{KenyonPemantle} (Theorem 7.8) by
an indirect method.

\begin{Theo}
  \label{theo:unic}
  For any formal initial condition $(g_x)_{x \in V(U)}$ where $U$ is a
  regular stepped solid, let $g_{(0,0,0)}$ be the solution at the
  origin of Kashaev's recurrence. Then:
  \begin{enumerate}
  \item $g_{(0,0,0)}$ is a Laurent polynomial in the $(g_x)_{x \in
      V(U)}$ vertex variables and in the face variables defined on the faces
    of $G(U)$;
  \item The monomials are in one-to-one correspondence with taut
    $C^{(1)}_2$ loop configurations on $U$;
  \item The $g$ variables appear with exponent in $\{-2, \dots,
    4\}$. The face variables appear with exponent in $\{0,1\}$;
  \item The coefficients in front of monomials are powers of $2$.
  \end{enumerate}
\end{Theo}

\begin{proof}
  For any $\sigma \in \Sigma(U)$ and $x\in V(U)$, the $g_x$ variable
  appears with an integer exponent in $w_U(\sigma)$. Indeed, it gets
  an exponent $\frac12$ when, and only when, the color of the loops
  change around $x$ (see the weights \eqref{eq:paramg} and
  \eqref{eq:weight}), and this happens an even number of times. The
  first point is thus a direct consequence of
  Theorem~\ref{theo:princ2}.
  
  The third point comes from the observation that any vertex belongs
  to at most $6$ faces of $G(U)$, so its exponent is between $0-2$ and
  $6-2$. The face variables can only appear once and with exponent
  $1$.

  The last point is a direct consequence of the second point, so all
  that remains to be proved is the following statement:

  Let $\sigma, \sigma' \in \Sigma(U)$ be two taut $C^{(1)}_2$ loop
  configurations on $G(U)$. If the following expressions in the formal
  $g$ variables are equal:
  \begin{equation}
    \label{eq:formws}
    \left( \prod_{(f,i_k)\in \sigma}w^f_{i} \right)
    \left( \prod_{x\in V(U)} \frac{1}{g_x^2}\right) \ = \
    \left( \prod_{(f,i_k)\in \sigma'}w^f_{i} \right)
    \left( \prod_{x\in V(U)} \frac{1}{g_x^2}\right)
  \end{equation}
  then $\sigma = \sigma'$.

  To prove this, we give a procedure to reconstruct $\sigma$ from the
  monomial in the left-hand side of \eqref{eq:formws}.

  Suppose that there is a vertex $x\in\Z^3$ so that $\sigma$ is
  already known on every face around $x$ except for one, which we call
  $f$. We claim that we can find $\sigma|_f$. To do so, first we look
  whether the face variable associated to $f$ is present in the
  monomial. If it is present, then the local configuration of $\sigma$
  at $f$ belongs to the third or fourth row of \eqref{eq:paramg};
  otherwise it belongs to the first, second or fifth row. Then we look
  at the exponent of $g_x$ that doesn't come from already known faces:
  in the first case it can be $\frac12$ (third row) or $0$ (fourth
  row). In the second case it can be $1$ (first row), $0$ (second row)
  or $\frac12$ (fifth row), so now we know which row $\sigma|_{f}$
  belongs to. There are two local configurations in this row. To know
  which one it is, just look at the color of an incoming edge that is
  already known.

  Since $\sigma$ is taut, we already know it outside of a sufficiently
  big ball centered at the origin. We want to use the previous
  argument to successively discover new faces, until $\sigma$ is known
  everywhere. To do so, we need to show that there is always an $x$
  that satisfies the first statement of the previous paragraph. Any
  $x$ having degree $2$ on the boundary of the graph formed by
  currently unknown faces would do the trick.

  We prove that such an $x$ always exists by showing a slightly more
  general result on \textit{lozenge graphs}. A lozenge graph is a
  finite planar quadrangulation such that all internal faces are
  non-degenerate rhombi with same edge length.  We call \emph{outer
    boundary} of a lozenge graph the set of edges separating a face
  from the infinite connected component of the complement of the graph
  in the plane.

  \begin{figure}[h]
    \centering
    \input{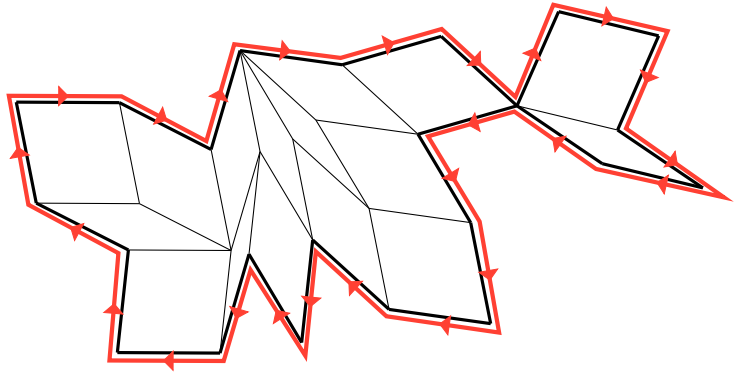}
    \caption{A connected lozenge graph and its oriented outer
      boundary.}
    \label{fig:exlos}
  \end{figure}

  \begin{Le}
    \label{lemma:deg2}
    Any non-empty lozenge graph has at least three
    vertices of degree $2$ on its outer boundary.
  \end{Le}

  \begin{proof}
    Let $G$ be a lozenge graph. By restricting ourselves to a
    connected component, we can assume that $G$ is connected. Let
    $x_1$ be a vertex on the boundary of $G$. We orient the boundary
    clockwise, meaning that we orient each of its edges so that the
    infinite connected component is on the left, see
    Figure~\ref{fig:exlos}. Starting from $x_1$, we follow the
    boundary by taking the leftmost edge when several outgoing edges
    are present. We denote $(x_1, \dots, x_p, x_1)$ this closed path,
    which is a single oriented curve around $G$, with possible pinches
    (some of the $x_i$ vertices may be equal). We also take
    $x_0 := x_p$.

    For any $i\in\{1,\dots,p\}$, the edge $\{x_{i-1},x_i\}$ belongs to
    exactly one rhombus (otherwise it would not be on the
    boundary). Let $u_i$ be the vertex of that lozenge that is
    diagonally opposite to $x_{i-1}$. Similarly, let $v_i$ be the
    vertex that is diagonally opposite to $x_{i+1}$ in the rhombus
    that contains the edge $\{x_i,x_{i+1}\}$; see
    Figure~\ref{fig:theta}.

    \begin{figure}[h]
      \centering
      \input{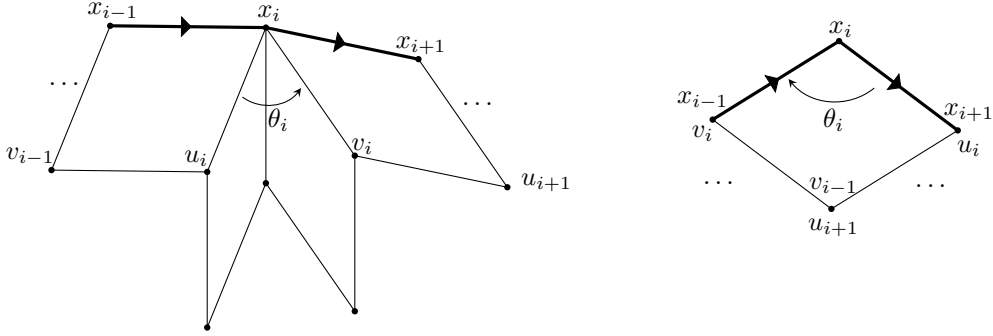}
      \caption{The $u_i$, $v_i$ vertices and $\theta_i$ angle when
        $x_i$ has degree $>2$ (left) or $2$ (right).}
      \label{fig:theta}
    \end{figure}

    Let $\theta_i$ be
    \begin{equation*}
      \theta_i = \angle (\overrightarrow{x_i x_{i-1}},
      \overrightarrow{x_i x_{i+1}}) -
      \angle (\overrightarrow{x_i x_{i-1}}, \overrightarrow{x_i
        u_{i}}) -
      \angle (\overrightarrow{x_i v_{i}}, \overrightarrow{x_i x_{i+1}}),
    \end{equation*}
    \emph{where the three angles are taken in} $]0,2\pi[$. Clearly,
    $\theta_i$ is equal to
    $ \angle (\overrightarrow{x_i u_{i}}, \overrightarrow{x_i v_{i}})$
    modulo $2\pi$. By checking all possible cases like in
    Figure~\ref{fig:theta}, we can be more precise:
    
    \begin{itemize}
    \item when $x_i$ does not have degree $2$ in $G$, $\theta_i \in
      ]0,2 \pi[$,
    \item otherwise, $\theta_i \in ]-\pi,0[$.
    \end{itemize}

    Our goal is to show that the sum of these angles is $-2\pi$, so
    that at least three of them have to be negative, and we can conclude.

    Notice that $\overrightarrow{x_iv_i} =
    \overrightarrow{x_{i+1}u_{i+1}}$. By using this fact and reorganizing
    the sums we get:
    \begin{equation*}
      \begin{split}
        \sum_{i=1}^p \theta_i & = \sum_{i=1}^p \angle
        (\overrightarrow{x_i x_{i-1}}, \overrightarrow{x_i x_{i+1}}) \
        - \sum_{i=1}^p \angle (\overrightarrow{x_i x_{i-1}},
        \overrightarrow{x_iu_{i}}) \ -
        \sum_{i=1}^p \angle (\overrightarrow{x_{i+1} u_{i+1}},
        \overrightarrow{x_i x_{i+1}})
        \\
        & = \sum_{i=1}^p \angle (\overrightarrow{x_i x_{i-1}},
        \overrightarrow{x_i x_{i+1}}) \ - \sum_{i=1}^p \left( \angle
          (\overrightarrow{x_i x_{i-1}}, \overrightarrow{x_i u_{i}})
          + \angle (\overrightarrow{x_{i} u_{i}},
          \overrightarrow{x_{i-1} x_{i}}) \right).
      \end{split}
    \end{equation*}

    By the choice of angles,
    $\angle (\overrightarrow{x_i x_{i-1}}, \overrightarrow{x_i u_{i}})
    + \angle (\overrightarrow{x_{i} u_{i}}, \overrightarrow{x_{i-1}
      x_{i}})$ is equal to $\pi$, so we have:
    \begin{equation*}
      \sum_{i=1}^p \theta_i =
      \sum_{i=1}^p \left( \angle (\overrightarrow{x_i x_{i-1}},
        \overrightarrow{x_i x_{i+1}}) - \pi \right).
    \end{equation*}

    Now the angles
    $\angle (\overrightarrow{x_i x_{i-1}}, \overrightarrow{x_i
      x_{i+1}}) - \pi$ are the oriented angles
    $\angle (\overrightarrow{x_{i-1}x_i}, \overrightarrow{x_i
      x_{i+1}})$ taken in $]-\pi,\pi[$. Since the boundary
    $(x_1, \dots, x_p, x_1)$ is a clockwise oriented closed curve,
    they sum up to $-2\pi$, so that
    \begin{equation*}
      \sum_{i=1}^p \theta_i = -2\pi
    \end{equation*}
    as claimed.
  \end{proof}
  This concludes the proof of Theorem~\ref{theo:unic}.
\end{proof}

\subsection{Limit shapes}

The existence of limit shapes is shown exactly as in
\cite{KenyonPemantle}. We just do the computations here to show how it
fits into our particular framework. See also
\cite{PetersenSpeyer,diFrancescoSotoGarrido,George} for similar proofs
in the case of the octahedron and cube recurrences.

For $v=(i,j,k) \in \Z^3$, we define its \textit{height} as $h(v) =
i+j+k$.

For $N\in \Z^+$, let $U_N$ be the following regular stepped solid:
\begin{equation*}
  U_N = \{ C_v  \mid v \in (\Z^-)^3, h(v) \leq -N \}.
\end{equation*}

\begin{figure}[!h]
  \centering
  \input{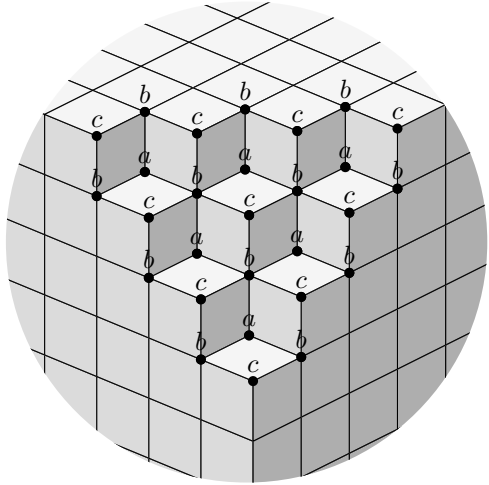}
  \caption{The stepped solid $U_3$ and periodic weights on the
    vertices.}
  \label{fig:un}
\end{figure}
We put a periodic $g$ function on the vertices $V(U_N)$ that only
depends on the height of vertices (see Figure~\ref{fig:un}). The
values of $g$ in the flat regions don't appear in the weight of a taut
configuration so they may be chosen arbitrarily. Our aim is to
describe the shape of a random taut loop configuration on $G(U_N)$
sampled proportionally to its weight, when $N$ is large.

\medskip

\begin{figure}[!h]
  \centering
  \input{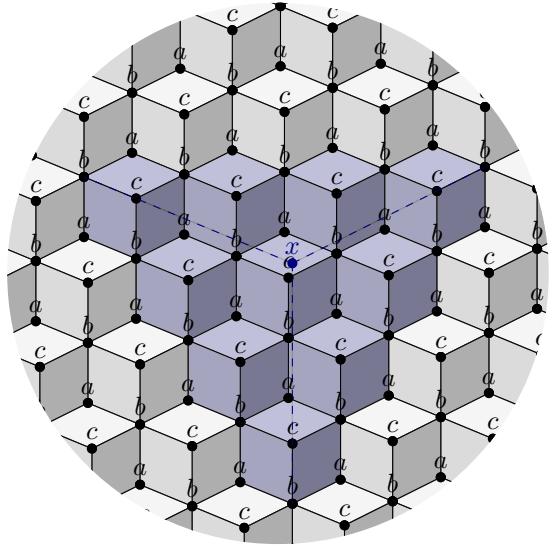}
  \caption{The infinite stepped solid $U$ with vertex weights
    $g^{a,b,c}$, an $x$ of height $5$ and the associated $U_x$ (in
    blue). The origin $(0,0,0)$ is one of the $a$ variables.}
  \label{fig:uinfi}
\end{figure}

Instead of letting $N$ change, it will be convenient to consider
instead the infinite stepped solid
\begin{equation*}
  U = \{ C_v  \mid h(v) \leq 2 \},
\end{equation*}
represented in Figure~\ref{fig:uinfi}, and to see it from some $x \in
\Z^3$ of positive height, and to let $x$ change. Thus for any $x \in
\Z^3$, we consider the ``regular'' stepped solid
\begin{equation*}
  U_x = U \cap (x+U_0).
\end{equation*}

Up to a translation of vector $-x$, $U_x$ is a regular stepped solid,
similar to $U_N$ where $N=h(x)-2$ for $h(x) \geq 3$.

Let $g^{a,b,c}$ denote the set of initial conditions of Figure
\ref{fig:uinfi}:
\begin{equation*}
  g^{a,b,c}(x) =
  \begin{cases} a \ \ \text{if} \ h(x) = 0, \\
    b \ \ \text{if} \ h(x) = 1, \\
    c \ \ \text{if} \ h(x) = 2.
  \end{cases}
\end{equation*}
Using these weights, we can define a partition function of loops on
$U_x$:

\begin{equation*}
  Y_x = \YY^{U_x}_{\text{taut}}(g^{a,b,c}).
\end{equation*}

Of course $Y_x$ depends only on the height of $x$. If $h(x)=N$, we
simply denote $Y_N = Y_x$. We also define
$X_N = \sqrt{Y_N Y_{N+2} + Y_{N+1}^2}$. These quantity can be exactly
computed using Kashaev's relation for a function depending only on
height, see \cite{KenyonPemantle}, Section 7.5. The result is the
following.

Let $R=\frac{ac}{b^2}$ and $S=\frac{bd}{c^2}$, where
$d=\frac{2b^3 + 3abc + 2(ac+b^2)^{\frac32}}{a^2}$. As a side note, $S$
and $R$ can be deduced from one another as the greatest root of the
intrinsic relation
\begin{equation*}
  R^2S^2 - 6RS - 4R -4S -3 = 0.
\end{equation*}

We have:
\begin{equation*}
  \begin{split}
    Y_{2n} & = a^{1-2n} b^{2n} R^{n^2} S^{n^2-n}, \\
    Y_{2n+1} & = a^{-2n} b^{2n+1} R^{n^2+n} S^{n^2}, \\
    X_{2n} & = \sqrt{1+R} \ Y_{2n+1}, \\
    X_{2n+1} & = \sqrt{1+S} \ Y_{2n+2}.
  \end{split}
\end{equation*}

We are interested in the quantity
\begin{equation}
  \label{defrho}
  \rho(x) = \left. \left( g_{(0,0,0)} \frac{\partial
        \log\left(\YY^{U_x}_{\text{taut}}(g)\right)}{\partial g_{(0,0,0)}}
    \right) \right|_{g=g^{a,b,c}}.
\end{equation}

This quantity has a probabilistic meaning. If $\sigma$ is a random
taut loop configuration on $G(U_x)$ chosen proportionally to its
weight $w^{U_x}_{\text{taut}}(\sigma)$ for the initial conditions
$g^{a,b,c}$, let $n_0$ be the power of $g_{(0,0,0)}$ appearing in the
formal weight of $\sigma$. For a face $f$, let $\epsilon_{f}(\sigma)$
be $1$ if $\sigma|_f$ is in the third or fourth row in
Figure~\ref{fig:locconf}, and $0$ otherwise. Then
\begin{equation}
  \label{eq:rhobserv}
  \rho(x) = \mathbb{E} \left[ n_0 +
    \frac{1}{2(1+R)} \sum_{f \in F \text{ around } (0,0,0)}
    \epsilon_f(\sigma) \right] .
\end{equation}
If we looked instead at $U_N$ for $N=h(x)-2$, by a simple symmetry,
this quantity would be equal to the same expectation on the vertex
$-x$ instead of $0$ (and this vertex would have to be of type $a$ in
$U_N$). By defining the same partial derivative as in \eqref{defrho}
with respect to $g_{(1,0,0)}$ or $g_{(1,1,0)}$, we could keep track of
similar observables for vertices of type $b$ and $c$. We will not make
use of the exact formula of that observable, we just use it to show
that the behavior of loops changes depending on the vertex of $U_N$.

The observable $\rho$ is defined as some logarithmic derivative of the
partition function $\YY^{U_x}_{\text{taut}}$. By taking the
logarithmic derivative of Kashaev's relation \eqref{eq:kasharel},
which is satisfied by $\YY^{U_x}_{\text{taut}}$, and evaluating at the
initial condition $g^{a,b,c}$, we get linear relations on $\rho$:
\begin{equation}
  \label{eq:line1}
  \begin{split}
      \text{if } h(x) \text{ is even}, \
      \rho(x+e_1+e_2+e_3) = &
      \alpha \rho(x) + \beta  \left( \rho(x+e_1) + \rho(x+e_2) +
        \rho(x+e_3)\right) \\
      & + \gamma \left( \rho(x+e_1+e_2) + \rho(x+e_1 + e_3) +
        \rho(x+e_2 +e_3)\right),
  \end{split}
\end{equation}

\begin{equation}
  \label{eq:line2}
  \begin{split}
    \text{if } h(x) \text{ is odd}, \ \rho(x+e_1+e_2+e_3) = & \alpha'
    \rho(x) + \beta' \left( \rho(x+e_1) + \rho(x+e_2) +
      \rho(x+e_3)\right) \\
    & + \gamma' \left( \rho(x+e_1+e_2) + \rho(x+e_1 + e_3) +
      \rho(x+e_2 +e_3)\right),
   \end{split}
\end{equation}
where
\begin{align*}
  \alpha &=\frac{3+3\sqrt{1+R}-2RS}{RS}, & \beta
  &=\frac{2+2\sqrt{1+R}+R}{R^2S}, & \gamma &=\frac{1+\sqrt{1+R}}{RS}, \\
  \alpha' &=\frac{3+3\sqrt{1+S}-2RS}{RS}, & \beta'
  &=\frac{2+2\sqrt{1+S}+R}{RS^2}, & \gamma' &=\frac{1+\sqrt{1+S}}{RS}.
\end{align*}

Let us define the generating function:
\begin{equation*}
  F(x,y,z) = \sum_{(i,j,k)\in \Z^3, h(i,j,k) \geq 0} \rho(i,j,k) x^i y^j z^k.
\end{equation*}

Using \eqref{eq:line1} \eqref{eq:line2}, it is straightforward to
compute $F$. It is a rational function of the form
\begin{equation*}
  F(x,y,z) = \frac{P(x,y,z)}{H(x,y,z)}
\end{equation*}
where $P$ is some polynomial and
\begin{equation*}
  \begin{split}
    H(x,y,z) = & \left(\alpha x y z + \gamma (x+y+z) \right)
    \left(\alpha' x y z + \gamma' (x+y+z) \right) \\
    & - \left(1 - \beta (xy + xz + yz) \right) \left(1 - \beta' (xy +
      yz + xz) \right).
  \end{split}
\end{equation*}

The coefficients $\alpha, \beta, \gamma, \alpha', \beta',
\gamma'$ can all be defined using
$R$ so they are all dependent. Actually, by defining
$\theta = \gamma \gamma'$, $H$ takes the form:

\begin{equation*}
  H(x,y,z) = \theta (x^2-1) (y^2-1) (z^2 -1) + (1-\theta) (xy-1)
  (xz-1) (yz-1).
\end{equation*}

At that point, we have recovered the denominator studied in
\cite{KenyonPemantle}. The asymptotic behavior of the observables
$\rho(i,j,k)$ can be obtained from the analysis of their generating
function $F$ at the singularity $(1,1,1)$ (Theorem 5.7 of
\cite{KenyonPemantle}, which is a corollary of various results of
\cite{PemantleWilson,BaryshnikovPemantle}). Around
that point, the leading homogeneous part of $H(1+X,1+Y,1+Z)$ is
\begin{equation*}
  \bar{H}(X,Y,Z) = 2(1+3\theta)XYZ + (1-\theta) (X^2Y + XY^2 + X^2Z +
  XZ^2 + Y^2Z + YZ^2).
\end{equation*}

Thus the limit shape of the model can be computed as the dual of the
curve
\begin{equation*}
  X^2Y + XY^2 + X^2Z + XZ^2 + Y^2Z + YZ^2 + \lambda XYZ
\end{equation*}
where $\lambda=\frac{2(1+3\theta)}{1-\theta}$. See
\cite{KenyonPemantle} for the computation of this dual, its shape, and
the behavior of $\rho$ depending on its position relatively to the
limit shape. The dual is a projective curve in $P\R^3$, in the
following figures we represent it in $\R^3$ intersected with the plane
$x+y+z=-1$.

\begin{figure}[h]
  \centering
  \includegraphics[width=7cm]{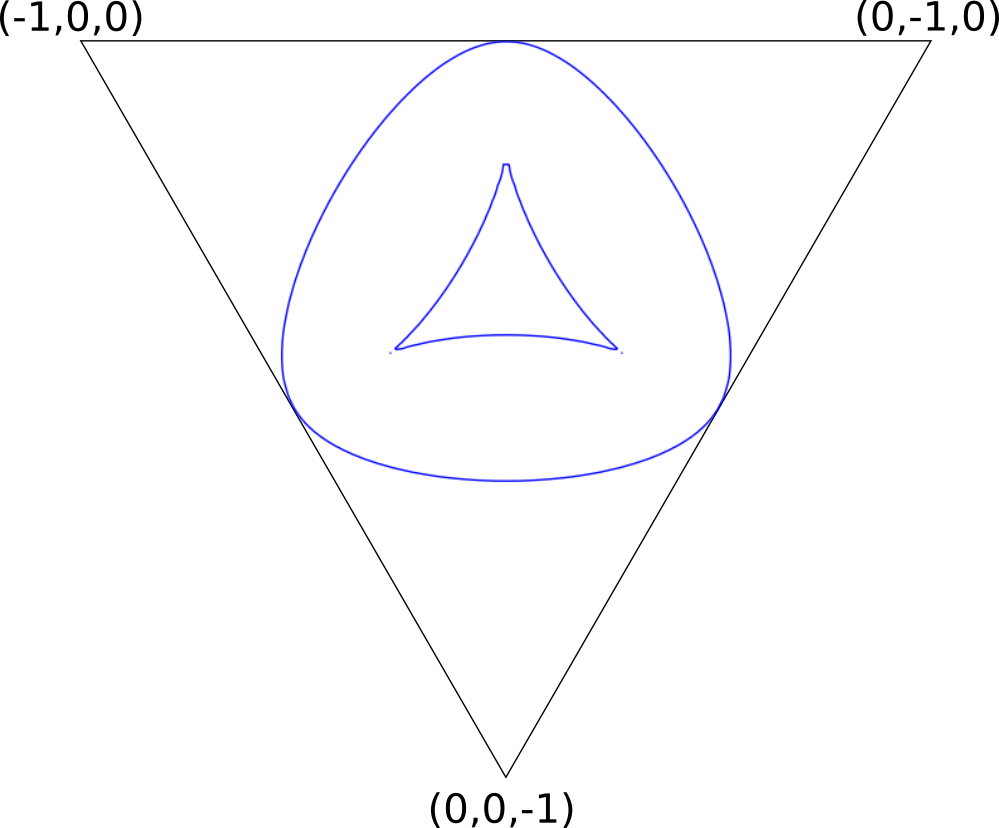}
  \hspace*{0.5cm}
  \includegraphics[width=7cm]{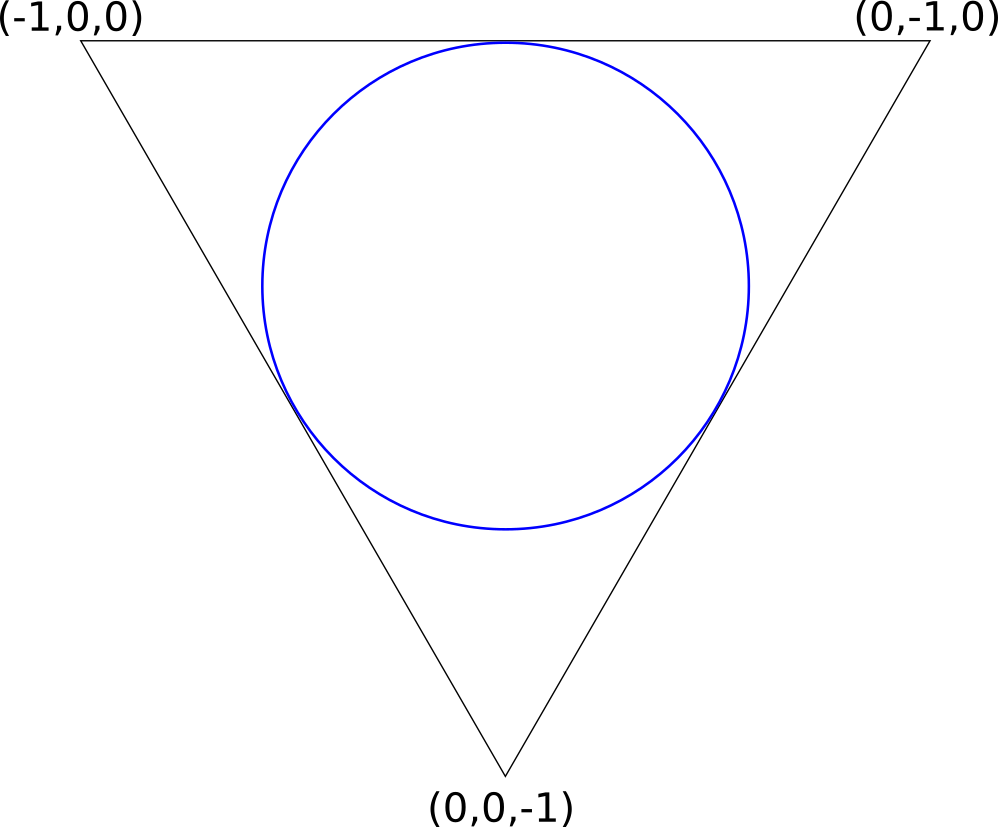}
  \caption{Limit shapes for $R=0.2$ (left) and $R=3$ (right).}
  \label{fig:limshape}
\end{figure}

We note that the limit shape only depends on $R$, and there is an
extra symmetry: when $R$ and $S$ are exchanged, $\lambda$ remains the
same so the limit shape is the same.  In general, $\lambda \in
(2,3]$. When $\lambda \neq 3$, the limit shape looks like a rounded
triangle tangent to the borders of the carved section of $U_N$, with
an internal facet; see Figure~\ref{fig:limshape}. The point
$\lambda=3$ is somehow critical and corresponds to $R=S=3$. At that
point, the limit shape becomes a circle and the central facet is
reduced to a point.

\begin{figure}[!h]
  \centering
  \includegraphics[height=7.5cm]{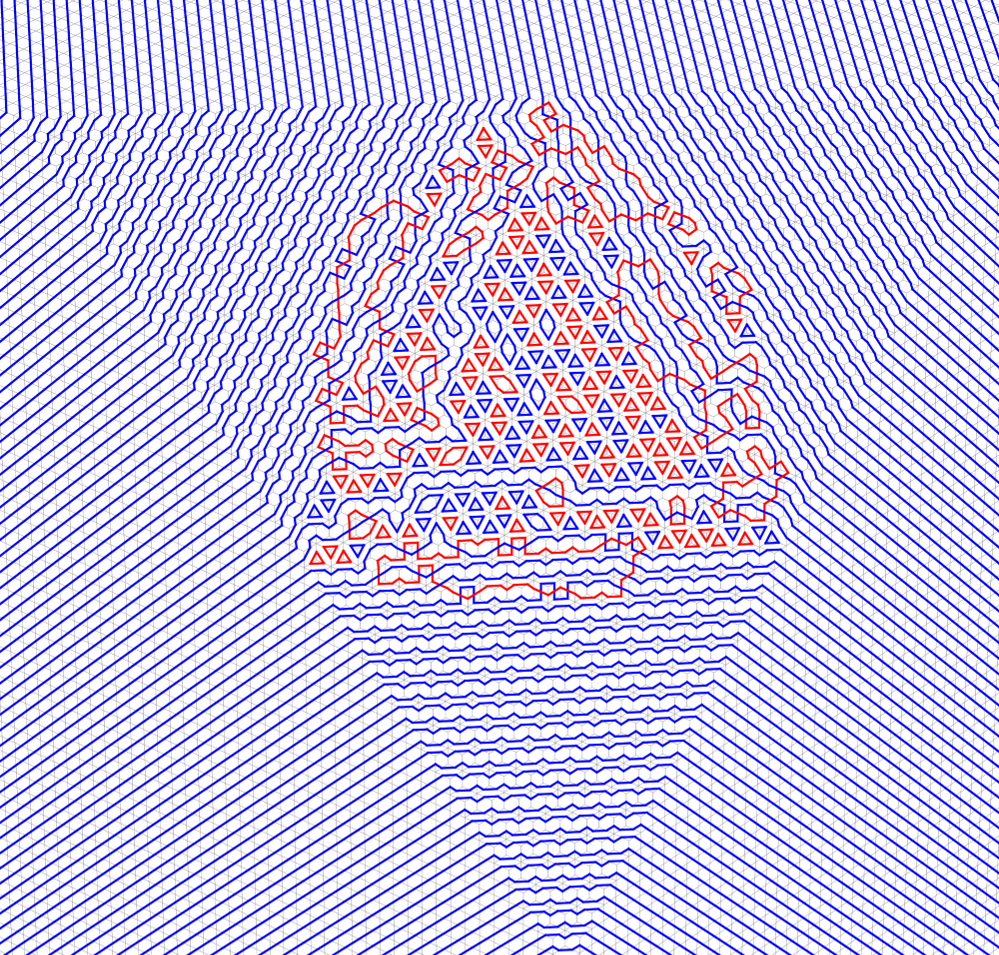}
  \includegraphics[height=7.5cm]{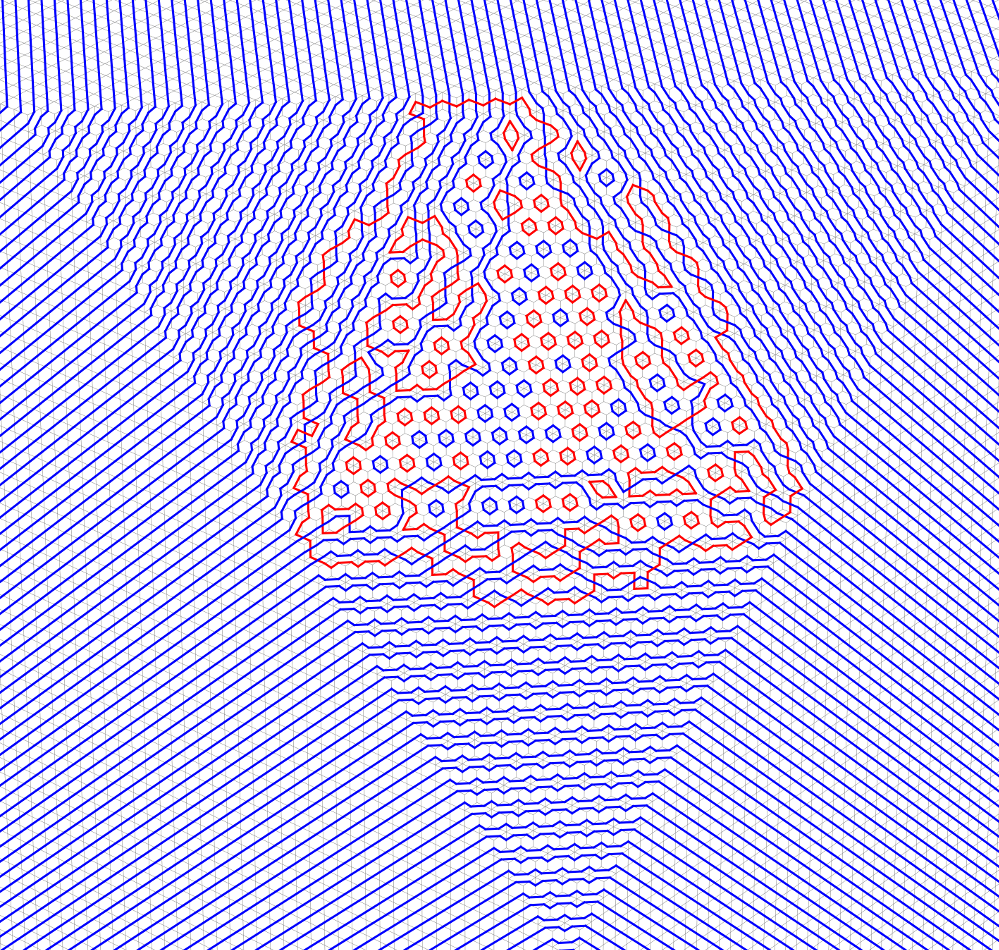}
  \includegraphics[height=7.5cm]{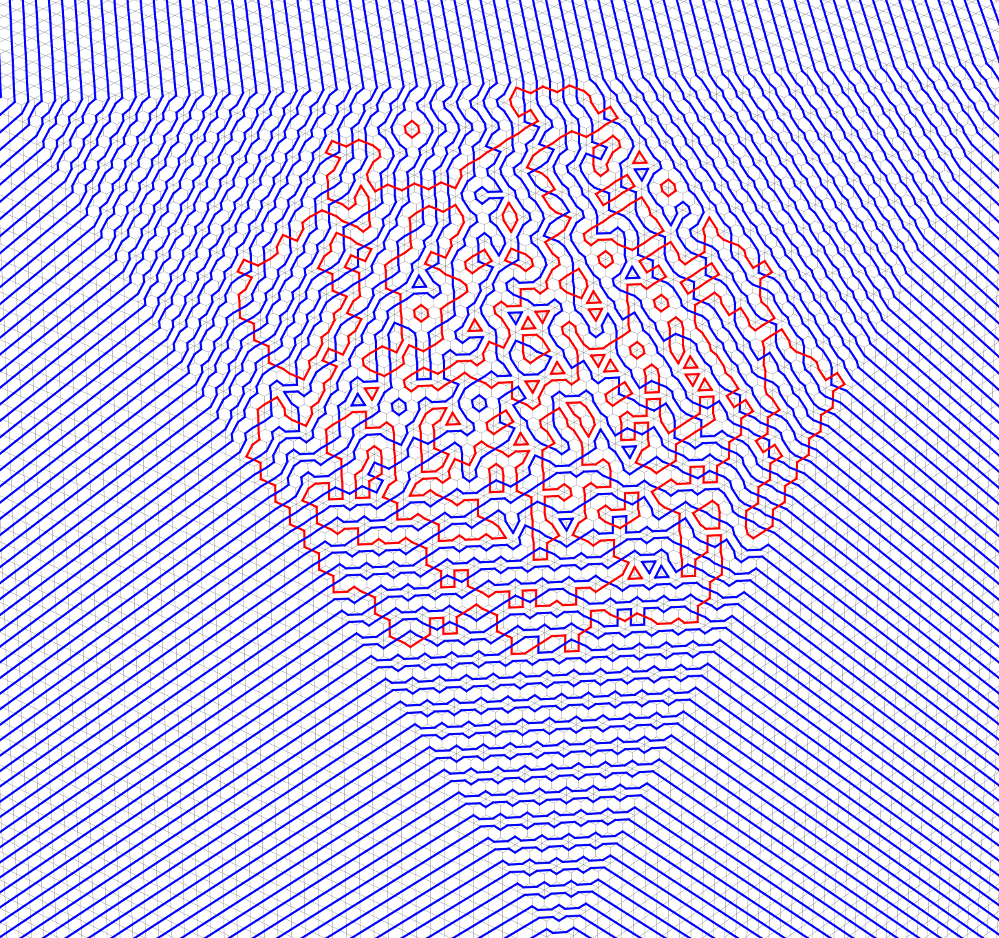}
  \caption{Simulations for $N=40$, and $R=0.2$ (top left), $R\simeq
    130.7$ ($S=0.2$) (top right), $R=3$ (bottom). The first two
    simulations correspond to the same limit shape.}
  \label{fig:simu}
\end{figure}

We computed simulations of the model for different values of $R$ in
Figure~\ref{fig:simu}. In the three corner regions $\rho$ decays to
$0$ exponentially fast in $N$, which corresponds to the ``frozen
phase'' where only infinite blue paths appear in the densest possible
packing; it is possible to convince oneself that this is indeed the
behavior implied by \eqref{eq:rhobserv} being close to $0$. The
annular region around the facet corresponds to a ``liquid phase''
where $\rho$ tends to $0$ polynomially fast. It seems that this
region's interface with the central facet is delimited by the infinite
blue paths closest to the center. In the central facet there is a
``gaseous phase'' where $\rho$ tends to $\frac13$, and the boundary
conditions don't appear any more.

\subsection{Cube groves}
\begin{figure}[h]
  \centering
  \input{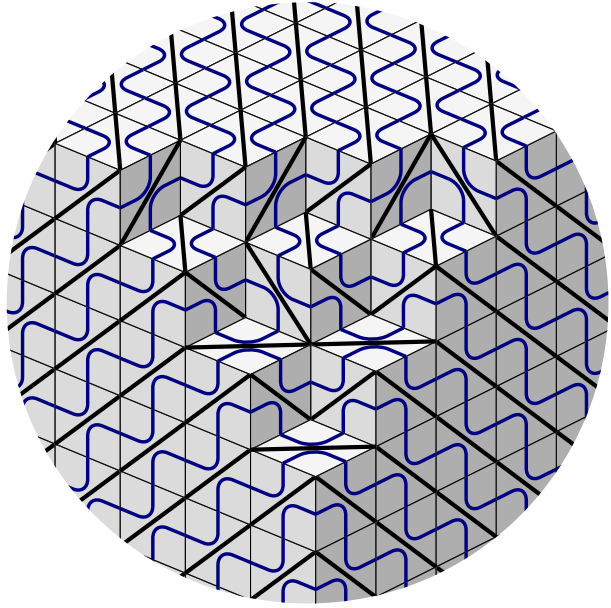}
  \caption{A configuration $\sigma \in \Sigma_0(U)$ and the
    corresponding cube grove.}
  \label{fig:grove}
\end{figure}

Cube groves were introduced by Carroll and Speyer in
\cite{CarrollSpeyer}.  They are essential spanning forests, often
represented with their dual forest, on the graph consisting of the
even vertices of $G(U)$ with edges on diagonals of the faces of the
cubes. One example is displayed in black lines in
Figure~\ref{fig:grove}.

\begin{figure}[h]
  \centering
  \input{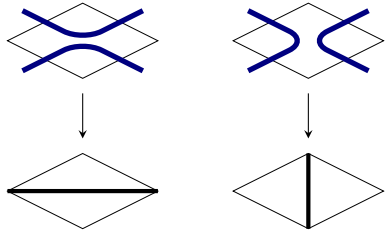}
  \caption{Local transformation from $\Sigma_0(U)$ to cube groves}
  \label{fig:loopsgroves}
\end{figure}

Let $\Sigma_0(U)$ be the subset of $\Sigma(U)$ containing all taut
$C^{(1)}_2$ loop configurations $\sigma \in \Sigma(U)$ such that
$N_{\sigma} = 0$, \textit{i.e.} $\sigma$ has no finite loop. Such a
configuration cannot contain any red edge, since all red paths are
finite loops. Thus $\sigma$ can be represented by a subset of the
diagonals of the faces of the cubes, as in
Figure~\ref{fig:loopsgroves}.

It is easy to check that this set of diagonals is necessarily a cube
grove, and conversely, any cube grove corresponds to such a loop
configuration (which is in fact the classical interface between the
spanning forest and its dual). The transformation is thus a bijection
between $\Sigma_0(U)$ and cubes groves on $U$. Moreover, this
bijection is weight-preserving, in the sense that the weight of
$\sigma\in \Sigma_0(U)$ is equal to the weight of its associated cube
grove as defined by Carroll and Speyer, when their face variables are
set to be equal to $1$.

As a result, the partition function truncated to $\Sigma_0(U),$
\begin{equation*}
  \YY^U_{0 \, \text{taut}}(g) = \sum_{\sigma \in \Sigma_0(U)} w_{\text{taut}}^U(\sigma),
\end{equation*}
is formally equal to the solution of the cube recurrence with initial
conditions on $U$. This is not such a surprise: on a field of
characteristic $2$, Kashaev's relation \eqref{eq:kasharel} reduces to
\begin{equation*}
  g g_{123} + g_{1} g_{23} + g_{2} g_{13} + g_{3} g_{12} = 0,
\end{equation*}
which is exactly the cube recurrence in characteristic $2$, while
$\YY^U_{\text{taut}}(g) $ reduces to $\YY^U_{0 \, \text{taut}}(g) $.

\appendix

\section{Calculations for Lemma~\ref{lemma:map}}
\label{app:map}

\begin{itemize}
\item $i=2$:
  \begin{equation*}
    \begin{split}
      \frac{1}{g_x^2}\sum_{\sigma \in L_2} w_{\loc}(\sigma) & =
      g_1 g_3 g_{12} g_{23},\\
      \frac{1}{(g'_x)^2}\sum_{\sigma \in L'_2} w_{\loc }(\sigma) & =
      g_1 g_3 g_{12} g_{23}.
    \end{split}
  \end{equation*}
  
\item $i=3$:
  \begin{equation*}
    \begin{split}
      \frac{1}{g_x^2}\sum_{\sigma \in L_3} w_{\loc}(\sigma) & =
      g^{-1} g_1 g_3 g_{12}^{\frac12} g_{23}^{\frac12} XZ \ +
      \ g^{-1} g_1 g_2 g_3 g_{12}^{\frac12} g_{23}^{\frac12} Y \\
      & = g_1 g_3 g_{12}^{\frac12} g_{23}^{\frac12} \left( \frac{XZ +
          g_2 Y}{g} \right),\\
      \frac{1}{(g'_x)^2}\sum_{\sigma \in L'_3} w_{\loc }(\sigma)
      & = g_1 g_3 g_{12}^{\frac12} g_{23}^{\frac12} Y_2.
    \end{split}
  \end{equation*}
  The equality of
  $\frac{1}{(g'_x)^2}\sum_{\sigma \in L'_2} w_{\loc }(\sigma)$ and
  $\frac{1}{g_x^2}\sum_{\sigma \in L_2} w_{\loc }(\sigma)$ is
  equivalent to Item \ref{fy} of Proposition \ref{prop:calculs}.
\item $i=4$:
  \begin{equation*}
    \begin{split}
      \frac{1}{g_x^2}\sum_{\sigma \in L_4} w_{\loc}(\sigma) & =
      g^{-1} g_1^{\frac12} g_2 g_3^{\frac12} g_{13}^{\frac12}
      g_{23}^{\frac12} YZ \ + \ g^{-1} g_1^{\frac32} g_2 g_3^{\frac12}
      g_{13}^{\frac12} g_{23}^{\frac12} X\\
      & = g_1^{\frac12} g_2 g_3^{\frac12} g_{13}^{\frac12}
      g_{23}^{\frac12} \left( \frac{YZ + g_1 X}{g} \right), \\
      \frac{1}{(g'_x)^2}\sum_{\sigma \in L'_4} w_{\loc}(\sigma) & =
      g_1^{\frac12} g_2 g_3^{\frac12} g_{13}^{\frac12}
      g_{23}^{\frac12} X_1.
    \end{split}
  \end{equation*}
  The equality is equivalent to Item \ref{fx} of Proposition
  \ref{prop:calculs}.
\item $i=5$:
  \begin{equation*}
    \begin{split}
      \frac{1}{g_x^2}\sum_{\sigma \in L_5} w_{\loc}(\sigma) & =
      2 g^{-2} g_1^{\frac32} g_2^2 g_3^{\frac32}Y \ + \ 2 g^{-2}
      g_1^{\frac32} g_2 g_3^{\frac32}XZ \ + \ g^{-1} g_1^{\frac12}
      g_2 g_3^{\frac12}g_{13}XZ \\
      & \ \ + g^{-1} g_1^{\frac32} g_2 g_3^{\frac12}g_{23}Y + g^{-1}
      g_1^{\frac12} g_2 g_3^{\frac32}g_{12}Y \\
      & = g_1^{\frac12} g_2 g_3^{\frac12} \left(\frac{2g_1g_2g_3 Y}{g^2} +
        \frac{2g_1g_3 XZ}{g^2} + \frac{g_{13}XZ}{g} +
        \frac{g_1g_{23}Y}{g} + \frac{g_3g_{12} Y}{g} \right),\\
      \frac{1}{(g'_x)^2} \sum_{\sigma \in L'_5} w_{\loc }(\sigma)
      & = g_1^{\frac12} g_2 g_3^{\frac12} X_1 Z_3.
    \end{split}
  \end{equation*}
  Using Items \ref{fx} and \ref{fz} of Proposition \ref{prop:calculs},
  we get
  \begin{equation*}
    \begin{split}
      X_1 Z_3 & = \left(\frac{g_1 X + Y Z}{g}\right)
      \left(\frac{g_3 Z + X Y}{g} \right) \\
      & = \frac{g_1 g_3 XZ + g_1 X^2
        Y + g_3 Y Z^2 + XY^2Z}{g^2} \\
      & = \frac{g_1g_3XZ + g_1 (g
        g_{23}+g_2g_3)Y + g_3 (g g_{12} + g_1 g_2) Y + (g g_{13} + g_1 g_3)
        XZ}{g^2}\\
      & = \frac{2g_1g_2g_3 Y}{g^2} + \frac{2g_1g_3 XZ}{g^2} +
      \frac{g_{13}XZ}{g} + \frac{g_1g_{23} Y}{g} +
      \frac{g_3g_{12}Y}{g}
    \end{split}
  \end{equation*}
  and the equality follows.
\item $i=6$:
  \begin{equation*}
    \begin{split}
      \frac{1}{g_x^2}\sum_{\sigma \in L_6} w_{\loc }(\sigma) & =
      g_1^{\frac12}g_2 g_3^{\frac12} g_{12}^{\frac12}g_{13}
      g_{23}^{\frac12},\\
      \frac{1}{(g'_x)^2}\sum_{\sigma \in L'_6}
      w_{\loc }(\sigma) & = g_1^{\frac12}g_2 g_3^{\frac12}
      g_{12}^{\frac12}g_{13} g_{23}^{\frac12}.
    \end{split}
  \end{equation*}
\item $i=7$:
  \begin{equation*}
    \begin{split}
      \frac{1}{g_x^2}\sum_{\sigma \in L_7} w_{\loc }(\sigma) & =
      g^{-1} g_1^{\frac12} g_3^{\frac12} g_{12}^{\frac12}
      g_{23}^{\frac12} XYZ \ + \ g^{-1} g_1^{\frac32} g_2
      g_3^{\frac32} g_{12}^{\frac12} g_{23}^{\frac12} \\
      & = g_1^{\frac12} g_3^{\frac12} g_{12}^{\frac12} g_{23}^{\frac12}
      \left(\frac{XYZ + g_{1} g_{2} g_{3} }{g} \right),\\
      \frac{1}{(g'_x)^2}\sum_{\sigma \in L'_7} w_{\loc }(\sigma) & =
      g_1^{\frac12} g_3^{\frac12} g_{12}^{\frac12} g_{23}^{\frac12}
      g_{123}^{-1} X_1 Y_2 Z_3 \ + \ g_1^{\frac12}
      g_3^{\frac12} g_{12}^{\frac32} g_{13} g_{23}^{\frac32} g_{123}^{-1} \\
      & = g_1^{\frac12} g_3^{\frac12} g_{12}^{\frac12}
      g_{23}^{\frac12}
      \left(\frac{X_1Y_2Z_3 + g_{12}g_{13}g_{23}}{g_{123}}
      \right).
    \end{split}
  \end{equation*}
  The equality is equivalent to Item \ref{ff} of Proposition
  \ref{prop:calculs}.
\end{itemize}

\section{Full Kashaev parametrization of free-fermionic $C^{(1)}_2$
  loops}
\label{app:fullc12}

Let $G$ be a planar quadrangulation with a boundary, $V$ its set of
vertices, $F$ its set of internal faces. Then $G$ is necessarily
bipartite; we fix a bipartite coloring of $V$ into black and white
vertices.

\begin{figure}[!h]
  \centering
  \def\svgwidth{5cm} \import{./}{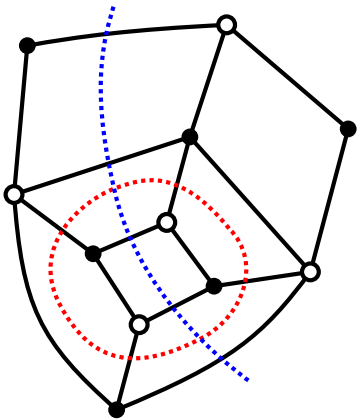}
  \caption{A planar finite quadrangulation with a boundary, and two
    train tracks.}
  \label{fig:tt}
\end{figure}

A \textit{train track} of $G$ is a path on the dual graph $G^*$
(defined with ``half-edges'' on the boundary of $G$, \textit{i.e.} not
connected at the external face), such that whenever it enters a face
it exits on the opposite edge of that face; see
Figure~\ref{fig:tt}. Let $T$ be the set of all train tracks of $G$.

It is a theorem of Kenyon and Schlenker \cite{KenyonSchlenker} that a
quadrangulation is a lozenge graph \textit{iff} no train track $t\in
T$ is a loop or crosses itself, and two distinct train tracks $t,
t'\in T$ cross at most once. For instance the quadrangulation of
Figure~\ref{fig:tt} cannot be made of non-degenerate rhombi.

Here we show that on a slightly more general class of
quadrangulations, an application is surjective. This linear
application will be needed to construct a parametrization such as
\eqref{eq:paramg}.

\begin{Le}
  \label{lemma:surj}
  If $G$ is a connected planar finite quadrangulation with a boundary,
  such that no train track is a loop, then the mapping
  \begin{alignat*}{3}
    \Phi : \ \ & \R^V && \rightarrow && \R^F \\
    &(h_v)_{v\in V} && \mapsto && \left( h_x + h_y - h_u -
      h_v\right)_{f \in F, \ f = \raisebox{-0.4cm}{\input{tikz/face.tex}}}
  \end{alignat*}
  is surjective.
\end{Le}

\begin{proof}
  Because of the rank-nullity theorem, it is sufficient to prove that
  the dimension of $\ker(\Phi)$ is $|V| - |F|$.
  
  Let $h$ be a vector in $\ker{\Phi}$, and let $t\in T$. We chose an
  orientation of $t$. Whenever $t$ crosses a face
  $ f = \raisebox{-0.4cm}{\input{tikz/face.tex}}$, for instance with $x,u$
  on its left and $y,v$ on its right, we have $h_x-h_v = h_u -
  h_y$. This means that the quantity given by the value of $h$ on the
  right minus its value on the left of an edge crossed by $t$ is
  constant along $t$. Let $\alpha_t(h)$ be this value.

  If we fix an orientation of every train track and a base vertex
  $x_0$, we thus get a linear transformation from $\ker(\Phi)$ to
  $\C^{|T|+1}$ by associating $(h_{x_0}, (\alpha_t(h))_{t\in T})$ to
  $h$. It is easy to see that this transformation is injective: if the
  family is null then $h_{x_0}=0$ and similarly for its neighbors
  using the train tracks adjacent to $x_0$, etc.
  
  Conversely, if we fix values $(\alpha_t)_{t\in T}$ and a $h_{x_0}$,
  we can to reconstruct a vector $h$ in the kernel of $\Phi$
  associated to these values: starting from $x_0$, define $h$ on its
  neighbors using the values $\alpha_t$ associated to train tracks
  adjacent to $x_0$, then on their neighbors, etc. The orientation of
  train tracks guaranties that these choices are coherent, and the
  fact that $\alpha_t$ is constant along $t$ is equivalent to
  $\Phi(h)$ being $0$ on every face crossed by $t$.

  Thus $\ker(\Phi)$ has dimension $|T|+1$.

  Let $E_{\text{ext}}$ be the set of edges adjacent to the external
  face, and $E_{\text{int}} = E \setminus E_{\text{ext}}$. Since the
  train tracks never loop, the set of train tracks define a coupling of
  external edges so
  \begin{equation}
    \label{eq:gr1}
    2 |T| = |E_{\text{ext}}|.
  \end{equation}
  Since every internal face is a quadrangle, we have
  \begin{equation}
    \label{eq:gr2}
    4 |F| = 2|E_{\text{int}}| + |E_{\text{ext}}|.
  \end{equation}
  Finally we have Euler's formula:
  \begin{equation}
    \label{eq:gr3}
    |V| - |E| + |F| = 1.
  \end{equation}

  Combining \eqref{eq:gr1}, \eqref{eq:gr2}, \eqref{eq:gr3} and
  $|E| = |E_{\text{int}}| + |E_{\text{ext}}|$ easily gives
  \begin{equation}
    |T| + 1 = |V| - |F|
  \end{equation}
  so $\dim(\ker \Phi) = |V| - |F|$ as needed.
\end{proof}

Now we can go back to the proof that parametrization \eqref{eq:paramg}
exists for any free-fermionic loop model on a graph $G$ that satisfies
the assumption of Lemma~\ref{lemma:surj}.

On every face $f$ we have a set of positive weights
$w^f_1,\dots, w^f_5$ that satisfy \eqref{eq:ff}. Let us define
$\kappa_f = w^f_5$ and
$R_f = \left( \frac{w^f_1}{\kappa_f} \right)^2$. We get:
\begin{equation}
  \label{eq:pdskappa}
  \begin{cases}
    w^f_1 = \kappa_f \sqrt{R_f} \\
    w^f_2 = \kappa_f \sqrt{\frac{1}{R_f}} \\
    w^f_3 = \kappa_f \sqrt{1+R_f} \\
    w^f_4 = \kappa_f \sqrt{1+ \frac{1}{R_f}}\\
    w^f_5 = \kappa_f.
  \end{cases}
\end{equation}

By Lemma~\ref{lemma:surj}, there is a function $h:V\rightarrow \R$
such that on every face $ f = \raisebox{-0.4cm}{\input{tikz/face.tex}}$,
\begin{equation*}
  \log(R_f) = h_x + h_y - h_u - h_v.
\end{equation*}

If we set $g_x = e^{h_x}$, we get
\begin{equation*}
  R_f = \frac{g_x g_y}{g_u g_v}
\end{equation*}
and \eqref{eq:pdskappa} becomes
\begin{equation}
  \begin{cases}
    w^f_1 = \kappa_f \sqrt{\frac{g_x g_y}{g_u g_v}} \\
    w^f_2 = \kappa_f \sqrt{ \frac{g_u g_v}{g_x g_y}} \\
    w^f_3 = \kappa_f \sqrt{ \frac{g_x g_y + g_u g_v}{g_u g_v}} \\
    w^f_4 = \kappa_f \sqrt{ \frac{g_x g_y + g_u g_v}{g_x g_y}}\\
    w^f_5 = \kappa_f.
  \end{cases}
\end{equation}

Multiplying all weights at a face by a same constant doesn't change
the relative weights of loop configurations. Here, if we multiply all
weights by $\frac{\sqrt{g_xg_yg_ug_v}}{\kappa_f}$, we get the
parametrization \eqref{eq:paramg}.

\begin{Rk}
  One can prove in exactly the same way that Kashaev's
  parametrization of the Ising model \cite{Kashaev} is possible whenever
  the underlying quadrangulation satisfies the assumption of
  Lemma~\ref{lemma:surj}. For instance, the Ising model on any finite
  isoradial graph admits a Kashaev parametrization.
\end{Rk}

\bibliographystyle{abbrv} \bibliography{main}

\end{document}